\documentclass[12pt]{article}

\newif\iffrench\frenchfalse
\usepackage[utf8]{inputenc}
\usepackage[T1]{fontenc}

\usepackage[dvipsnames]{xcolor}

\iffrench
\usepackage[french]{babel}
\else
\usepackage[english]{babel}
\fi

\usepackage{pslatex}
\usepackage{lmodern}

\usepackage{amssymb}
\usepackage{amsthm}
\usepackage{bbold} 
\usepackage{dsfont}
\usepackage{mathtools}
\usepackage{thmtools}
\usepackage{soul,bbm}
\usepackage{paralist}

\usepackage{varioref} 
\usepackage{hyperref} 
\hypersetup{
    colorlinks=true,
    linkcolor=blue,
    }

\usepackage{tikz}
\usepackage{pstricks-add}

\usepackage[babel=true]{csquotes} 
\usepackage{setspace}
\usepackage{array}
\usepackage{layout}

\usepackage{stmaryrd} 
\usepackage{mathrsfs} 

\usepackage[margin=20pt,font=small,labelfont=bf]{caption} 

\usepackage[backend=bibtex,
  citestyle=numeric-comp,
  giveninits=true,
  url=false,
  isbn=false,
  doi=false,
  eprint=false,
  maxbibnames=50,
]{biblatex}
\newcommand{\smarge}[2]{\usepackage[top=#1,bottom=#1+1cm,left=#1-#2,right=#1]{geometry}}

\iffrench

\addto\captionsfrench{} 

\newcommand{\ioo}[1]{\left]#1\right[}
\newcommand{\ifo}[1]{\left[#1\right[}

\newtheorem{thm}{Theorème}[section]
\newtheorem{ppn}[thm]{Proposition}
\newtheorem{cor}[thm]{Corollaire}
\newtheorem{lem}[thm]{Lemme}
\newtheorem{dfi}[thm]{Définition}
\newtheorem{cjt}[thm]{Conjecture}

\newcommand{\pth}[1]{\left(#1\right)}
\newcommand{\cro}[1]{\left[#1\right]}
\newcommand{\acc}[1]{\left\{#1\right\}}
\newcommand{\abs}[1]{\left|#1\right|}
\newcommand{\dabs}[1]{\left\|#1\right\|}

\else

\newtheorem{thm}{Theorem}[section]
\newtheorem{ppn}[thm]{Proposition}

\newtheorem{lem}[thm]{Lemma}
\newtheorem{dfi}[thm]{Definition}

\theoremstyle{remark}
\newtheorem{req}[thm]{Remark}

\newcommand{\ioo}[1]{\left(#1\right)}
\newcommand{\ifo}[1]{\left[#1\right)}

\newcommand{\pth}[1]{\left(#1\right)}
\newcommand{\cro}[1]{\left[#1\right]}
\newcommand{\acc}[1]{\left\{#1\right\}}
\newcommand{\abs}[1]{\left|#1\right|}
\newcommand{\dabs}[1]{\left\|#1\right\|}

\newcommand{\eg}{e.g.~}

\fi

\newcommand{\esp}{\hspace{1cm}}

\def\BS{\color{Bittersweet}}

\newcommand{\tq}{\hspace{0.25cm}/ \hspace{0.25cm}}

\newcommand{\goq}{\geqslant}
\newcommand{\loq}{\leqslant}
\newcommand{\eps}{\varepsilon}

\newcommand{\de}{\,\mathrm{d}}

\newcommand{\dr}{\partial}

\DeclareMathOperator\arcoth{arcoth}

\newcommand{\Er}{\mathds{R}}
\newcommand{\Zed}{\mathds{Z}}
\newcommand{\N}{\mathds{N}}

\newcommand{\mcc}{\mathcal{C}}

\newcommand{\mcl}{\mathcal{L}}

\newcommand{\mct}{\mathcal{T}}

\smarge{2.4cm}{0cm}

\usepackage{enumitem}
\usepackage{soul}

\setstcolor{red}

\renewbibmacro{in:}{}
\renewcommand{\leq}{\leqslant}
\renewcommand{\geq}{\geqslant}

\title{Reaction-diffusion model for a population structured in 
 phenotype and space\\ I -- Criterion for persistence}

\author{Nathanaël Boutillon$^{\text{a,b}}$ and Luca Rossi$^{\text{c}}$ \\
\footnotesize{$^{\text{a}}$ Aix Marseille Univ, CNRS, I2M, Marseille, France}\\
\footnotesize{$^{\text{b}}$ INRAE, BioSP, 84914, Avignon, France}\\
\footnotesize{$^{\text{c}}$ Istituto \enquote{G. Castelnuovo}, Sapienza Università di Roma, Rome, Italy}
}

\begin{document}

\maketitle

\begin{abstract}
  We consider a reaction-diffusion model for a population structured in 
  phenotype. We assume that the population lives in a  
  heterogeneous periodic environment, so that a given phenotypic trait 
  may be more or less fit according to the spatial location. 
  The model features spatial mobility of individuals as well as mutation.
  
  We first prove the well-posedness of the model. 
  Next, we derive a criterion for the persistence of the population
  which involves the 
  {\em generalised principal eigenvalue} associated with 
  the linearised elliptic operator. 
  This notion allows us to handle the possible lack of coercivity 
of the operator.
  We then obtain a monotonicity result for the generalised principal eigenvalue, 
  in terms of the frequency of spatial
  fluctuations of the environment and in terms of the spatial diffusivity. 
  We deduce that
  the more heterogeneous is the environment, 
  or the higher is the mobility of individuals,
  the harder is the persistence for the species. 
  
  This work lays the mathematical foundation to investigate some 
  other optimisation problems
  for the environment to make persistence as hard or as easy as possible,
  which will be addressed in the forthcoming companion paper.
\end{abstract}



 \tableofcontents


\vspace{2cm}
 
\section{Introduction}

\subsection{Presentation of the model}

 We consider a model of population dynamics which takes into account 
both the spatial heterogeneity of the environment and the phenotypic structure of the population. This amounts to a reaction-diffusion equation with
a nonlocal reaction term of the Fisher-KPP type.

More precisely, let $u(t,x,\theta)$ denote the density at time $t$ of a population structured in space ($x$-variable) and in phenotype
($\theta$-variable). 
The domain of the spatial variable is the whole space $\Er^N$, $N\in\N$, and
the domain of the phenotype, denoted by~$\Theta$, is either $\Er^P$, $P\in\N$, 
or a smooth bounded domain of $\Er^P$. In the latter case we call $\nu$ the outward normal to~$\dr\Theta$.

The dynamic of the density is governed by the system
\begin{equation}\label{eq:main}
  \left\{
  \begin{aligned}
    \dr_tu&=\Delta_xu+\Delta_\theta u+u\pth{r(x,\theta)-\rho(t,x)},&t> 0,\  (x,\theta)\in\Er^N\times\Theta,\\
    \nu\cdot\nabla_{\theta} u&=0,&t> 0,\  (x,\theta)\in\Er^N\times\dr\Theta,\\
  \end{aligned}
  \right.
\end{equation}
with
\begin{equation}\label{eq:rho}
\rho(t,x):=\int_{\Theta} u(t,x,\sigma)\de \sigma
\end{equation}
representing the total population at given time and location.
We supplement problem \eqref{eq:main}-\eqref{eq:rho} with an initial condition
\begin{equation}\label{eq:IC}
u(0,x,\theta) = u_0(x,\theta), \qquad(x,\theta)\in\Er^N\times\Theta.
\end{equation}
The quantity $r(x,\theta)$ corresponds to the fitness of an individual of phenotype $\theta$ at location $x$, in the absence of competition. 
We refer to~$r$
as the \emph{space/phenotype fitness landscape}. We will assume 
in the sequel that $r$ is periodic in the spatial variable~$x$.
The term $\rho(t,x)u(t,x,\theta)$, which is nonlocal in the $\theta$ variable,
is the competition faced by individuals of phenotype $\theta$
located at $x$.
Throughout the paper, any reference to $\partial\Theta$, such as 
the boundary condition in~\eqref{eq:main}, has to be ignored in the
case $\Theta=\Er^P$.

The first question that we address in this work is the well-posedness of the problem,
which is not covered by the existing literature under the generality 
of our assumptions. 
 We will derive it by first showing that solutions remain bounded,
despite the equation lacks the comparison principle.

The next two fundamental questions on the model are:
\begin{enumerate}
\item {\em When does the population survive or get extinct?}
\item {\em How to optimise the chances of persistence/extinction?}
\end{enumerate}
In Section \ref{ss:intro_ltb}, we shall answer the first question 
by giving a persistence
criterion based on the {\em generalised principal eigenvalue} 
of the linearised operator.
Using this criterion, we shall provide 
in Section \ref{ss:intro_period} an answer to the second question
about the chances of persistence in terms of the fluctuations of the 
environment  and of the diffusivity of the species.
These results will be obtained by considering a fitness function 
$r(x/L,\theta)$ and by adding a diffusion coefficients $d$ in front of $\Delta_x$,
then modulating the parameters $L,d$. In loose terms, 
we will show that 
persistence becomes harder when the heterogeneity of the environment increases,
as well as when the diffusivity grows.
Other optimisation problems, under some specific structural assumptions on the 
landscape, will be addressed in the forthcoming companion work.

\subsection{Related models}
Let us explain where the model~\eqref{eq:main} comes from, and relate it to the existing
literature.

\paragraph{The classical Fischer-KPP equation.} 
The following equation has been independently introduced 
 in 1937 by Fisher \cite{Fis37} and Kolmogorov, Petrovsky and Piskunov \cite{KPP37}\,:
\begin{equation}\label{eq:kpp_classique}
  \dr_tu(t,x)=\Delta_xu+ru-u^2,\esp t>0,\esp x\in\Er.
\end{equation}
Since then, it has been widely studied, both from a mathematical point of view (\eg \cite{AroWei75}) and for applications in population dynamics (\eg \cite{Ske51}). 
Let us imagine that $u$ is the distribution of a population living on $\Er$. The term $\Delta_x u$ then stands for the movements of the individuals, the term $ru$ stands for the demography in the absence of competition, 
and the nonlinear term $-u^2$ encodes the competition pressure.

 When the individuals of a population move and reproduce, it can be expected that they will eventually invade their environment.
One of the most important aspects of the Fisher-KPP equation is that it shows an invasion taking place 
at the asymptotic linear speed $c_{KPP}:=2\sqrt r$, in the following sense.
If the initial condition is bounded, nonnegative, not identically equal to zero and has a compact support,
then: 
\[
\forall c\in\ioo{0,c_{KPP}},\quad
\lim_{t\to+\infty}u(t,ct)=r,\]
and 
\[
\forall c>c_{KPP},\quad
\lim_{t\to+\infty}u(t,c t)=0.\]
The standard Fisher-KPP equation has been extended in several ways, in order to make it more realistic and
adapted to capture some relevant biological properties. 
Our model~\eqref{eq:main} tries to circumvent two of the limitations of the standard equation.
First, the environment in which individuals live is rarely homogeneous \cite{ShiKaw97,SKT86,OkuLev01}. This is why the environmental conditions
should affect the growth rate~$r$, according to the spatial position~$x$. Second, the population itself is heterogeneous: for example, individuals
with distinct phenotypes may react differently to the same environmental conditions. The interplay between the different phenotypes can have a true impact on the dynamics of the population \cite{EllCor12,MBC19}.

In model~\eqref{eq:main}, the phenotype is taken into account through the new variable $\theta$.
Then, all
individuals located at the same position compete with each other, 
regardless of their phenotypes. 
Thus, the nonlinear term $-u^2$ is replaced by the nonlocal term $-\rho u$. 
Finally, we add a term $\Delta_{\theta}u$, to account for
mutations (which induce variations in phenotype); this modelling means that the mutations are small and frequent.

Models similar to~\eqref{eq:main} appeared for the first time in the work of Prévost~\cite{Pre04}. Champagnat and Méléard~\cite{ChaMel07} derived variants of~\eqref{eq:main} as large population limits of individual-based models (in their models, the mutations are rarer but larger than in ours). Their derivations allow one to understand better where the different terms of the equation come from. 

We now go back to the standard Fisher-KPP equation, and we progressively add features to it until we get equations similar to~\eqref{eq:main}, in order to understand the different mathematical properties that these features entail.

\paragraph{Heterogeneous landscape.}

If the demography is affected by the environment, and changes from place to place,
then one is led to consider a growth rate that depends on the spatial variable, i.e.\
$r=r(x)$.
This results in the following $N$-dimensional 
spatially heterogeneous version of the Fisher-KPP equation (typically $N=2$ in applications):
\begin{equation}\label{eq:kpp_hetero}
  \dr_tu(t,x)=\Delta_x u+r(x)u-u^2,\esp t>0,\esp x\in\Er^N.
\end{equation}
We have seen that in the homogeneous setting~\eqref{eq:kpp_classique},
provided $r>0$ (i.e.~the environment is favourable to the species) and the initial condition is nonzero and nonnegative,
the population persists and actually
invades the whole environment. This is called the \emph{hair-trigger effect}.

The model~\eqref{eq:kpp_hetero} allows one to consider environments composed also of unfavourable
regions, that is, $r(x)$ can be negative for some~$x\in\Er$.
Then, the question of the \emph{persistence} of the population is no longer trivial.
A sharp criterion to determine whether the population persists or, on the contrary,
eventually gets extinct, is provided by the sign 
of the principal eigenvalue of the operator linearised around the
steady state $u\equiv0$, i.e. the operator
\[\Delta_x+r(x).\]
(More details about principal eigenvalues are given in Section~\ref{ss:intro_ltb}
below.)
This criterion was first derived in~\cite{CanCos89}
for problems set in bounded domains. When dealing with unbounded environments,
a typical modelling assumption is that the landscape is {\em periodic}, see e.g.~\cite{SKT86}.
This can be viewed as a mathematical approximation of more complex situations.
In such a framework, the characterisation in terms of the principal eigenvalue
of an operator
acting on the space of periodic functions has been derived in ~\cite{BHR05}.

As in the homogeneous setting, for the heterogeneous equation~\eqref{eq:kpp_hetero} with periodic $r$,
when there is persistence there is also complete invasion towards a positive steady state.
The invasion occurs with some spreading speed $c_{het}$, which replaces the quantity $c_{KPP}$. 
When $N\geq 2$, the spreading speed $c_{het}$ may depend on the direction of spreading.
A formula for $c_{het}$ has been derived by Gärtner and Freidlin~\cite{GarFre79} 
using probabilistic techniques, and later with PDEs techniques, along with some extensions, see e.g.~\cite{Wei02,BHN08,Ros17}.
A general survey of mathematical works and methods for reaction-diffusion equation in heterogeneous environments is given~in~\cite{Xin00}.

\paragraph{Structured population.}
 Consider now a heterogeneous population structured according to
 the phenotype. A first situation is when there is a finite number of phenotypes: this leads to a system of Fisher-KPP equations, each of them accounting for
 a given phenotype, competing with the other ones. In this direction, let us mention \cite{Gir17} for homogeneous environments and \cite{Gir23,GriMat21} for heterogeneous environments.

In the present paper, we shall use another formalism, and assume that the phenotype is a continuous parameter. Let us first stick to a homogeneous environment; we get a simplified version of~\eqref{eq:main}:
\begin{equation}\label{eq:main_homo_intro}
  \left\{
  \begin{aligned}
    \dr_tu&=\Delta_xu+\Delta_{\theta} u+u\pth{r(\theta)-\rho(t,x)},&t>0,\  (x,\theta)\in\Er^N\times\Theta,\\
    \nu\cdot\nabla_{\theta} u&=0,&t>0,\  (x,\theta)\in\Er^N\times\dr\Theta,\\
    u(0,x,\theta)&=u_0(x,\theta),&(x,\theta)\in\Er^N\times\Theta,
  \end{aligned}
  \right.
\end{equation}
with $\rho$ defined by \eqref{eq:rho}.
Here, contrarily to~\eqref{eq:main}, the coefficient $r$ depends on $\theta$ only.  
This model is studied in \cite{BJS16} in the case $\Theta=\Er^P$, 
under the assumption that $-r$ is
coercive. The latter condition allows the authors to define in a standard way 
the principal eigenvalue $\lambda$ of the~operator
\[\Delta_x+\Delta_{\theta}+r(\theta)\]
acting on a Hilbert space. 
One then has a dichotomy analogous to that for the spatially heterogeneous 
Fisher-KPP equation
\eqref{eq:kpp_hetero}: if $\lambda<0$, then the population eventually gets extinct,
while if $\lambda>0$, then the population persists.
Moreover, in the latter case,
the homogeneity of the environment (i.e. the independence in $x$ of $r$) and the coercivity of $-r$ 
entail that the equation admits a unique positive steady state $V(\theta)$,
despite the lack of comparison principle, and  any nontrivial solution 
converges towards $V$ with a spreading speed equal to
$c_{phen}:=2\sqrt{\lambda}$.

\paragraph{Recent models.} Several spatially-heterogeneous versions of~\eqref{eq:main_homo_intro} have been studied
in the literature. In \cite{ACR13,Pel20}, the authors consider
the case of a population evolving along an environmental gradient, that is,
\[r(x,\theta)=R(x-B\theta),\esp B>0,\]
where the function $R:\Er\to\Er$ is positive on the interval $(-1,1)$, negative outside $[-1,1]$, and~$R(y)$
tends quadratically to $-\infty$ as $|y|\to\infty$.
In this model, the fittest phenotype at position $x$ lies on the interval $((x-1)/B,(x+1)/B)$, so the phenotype is driven by the spatial variable~$x$.
As a consequence, evolution is necessary for the spatial
propagation of the population. This is relevant when one studies, for example,  migrations along the north-south axis, or uphill-downhill for populations living in the mountains. In \cite{ACR13} again, using the coercivity
of $-r$, the authors give a criterion based on a principal eigenvalue to determine whether the population persists or gets extinct.

Additionally to spatial heterogeneity, the authors of \cite{ABR17} add a heterogeneity in time corresponding to climate change. Their growth rate takes the form $r(x-ct,\theta)$.
Different shapes of the favourable space/phenotype zone are studied (confined zone, environmental gradient, or a mixing of the two). In each case, the authors find a criterion for persistence and propagation of the species. The criteria, again, are based on the principal eigenvalue of linear operators.

There seems to be few works considering the phenotype together with a {periodic} heterogeneous environment. In \cite{AlfPel22}, Alfaro and Peltier focus on homogeneous environments with a very small periodic perturbation. When the perturbation is small enough, they show the existence of stationary states and pulsating travelling waves. In \cite{BouMir15}, Bouin and Mirrahimi focus on a situation where evolutionary events (mutations) occur on a much faster timescale than ecological events (movements).
In a similar vein, but in the stationary setting, Léculier and Mirrahimi~\cite{LecMir23} studied the weak mutation limit of the stationary version of a model related to ours. They also give biological properties of the solution; in some cases, they are able to determine if the population at equilibrium is polymorphic or monomorphic.

%
%
%

\section{Main results}


\subsection{Assumptions on the landscape}

The set $\Theta$ is either $\Er^P$ or a ${C}^{2}$ bounded domain in 
$\Er^P$. 
We shall make the following assumptions on $r$:
\begin{equation}\label{eq:general_assumptions}
  \left\{
  \begin{aligned}
    &r\in L^{\infty}_{loc}(\Er^N\times{ \overline\Theta}),\text{ $r$ is bounded above;}\\
    &r(x+h,\theta)=r(x,\theta)\text{ \;for all $(x,\theta)\in\Er^N\times\Theta$ and $h\in \Zed^N$.}\\
  \end{aligned}
  \right.
\end{equation}
Throughout the paper, for a given function defined on $\Er^N\times\Theta$,
we shall refer to the second condition 
in~\eqref{eq:general_assumptions} as ``$1$-periodicity in $x$''.
We point out that in the case where $\Theta$ is bounded, the function
$r$ is globally bounded.
In the case $\Theta=\Er^P$, we will sometimes assume that
    \begin{equation}\label{eq:r<0}
    \limsup_{\|\theta\|\to\infty}\;r(x,\theta)\loq 0,\quad
    \text{uniformly with respect to $x\in\Er^N$,}
    \end{equation}
    which means 
    $$\lim_{R\to+\infty}\left(\sup_{x\in\Er^N,\ \|\theta\|\goq R}r(x,\theta)\right)
    \loq0.$$
We point out that, unlike in other mathematical studies in the literature, 
we do not assume~$-r$ to be coercive. 
This will entail a loss of compactness
that we will face using the notion of {\em generalised principal eigenvalue}.

 \subsection{Well-posedness of the problem}\label{sec:well}
   We now state the well-posedness result for \eqref{eq:main}-\eqref{eq:IC}.
Solutions are understood in the weak sense, namely: $u\in L^{\infty}(\ifo{0,T}\times\Er^N\times\overline\Theta)\cap L^{\infty}(\ifo{0,T}\times\Er^N;L^1(\Theta))$ for all $T>0$, and
  it holds
  \begin{multline}\label{eq:weak_formulation}
    \forall\phi\in{C}^{\infty}_c(\ifo{0,+\infty}\times\Er^N\times\overline\Theta),\\
      -\int u(\dr_t\phi)-\int_{\Er^N\times\Theta} u_0(\cdot,\cdot)\phi(0,\cdot,\cdot)=\int u\pth{\Delta_x\phi+\Delta_{\theta}\phi}+\int u\phi\pth{r-\int_{\Theta}u(t,x,\cdot)}
  \end{multline}
  (the non-specified integrals are taken on $\ifo{0,+\infty}\times\Er^N\times\Theta$).
In fact, solutions to the Cauchy problem \eqref{eq:main}-\eqref{eq:IC}
satisfy the first equation in \eqref{eq:main}
also in the strong sense, i.e., they belong to
$W^{1,2}_{p,loc}(\ioo{0,+\infty}\times\Er^N\times\overline\Theta)$
for all $p>1$,
where, here and in the sequel, $W^{1,2}_p$ (resp.~$W^{1,2}_{p,loc}$) is the set of functions that belong to $L^p$ (resp. $L^p_{loc}$)
together with their first order derivative with respect to $t$
and their partial derivatives up to order $2$ with respect to $(x,\theta)$.
 In particular,  solutions 
 belong to $W^{1}_{p,loc}(\ioo{0,+\infty}\times\Er^N\times\overline\Theta)$
with respect to all variables $(t,x,\theta)$.
Hence, by Morrey's inequality, the solutions
are continuous, 
  so \eqref{eq:rho} holds for every $(t,x)\in \ifo{0,+\infty}\times\Er^N$
  (and not just almost everywhere).
  Moreover, since the initial datum will be assumed to be continuous, the initial 
  condition is attained in the classical sense too.
  However, without further regularity assumptions on $r$,
  the solution $u$ cannot be expected to be more than Lipschitz continuous in the $(x,\theta)$ 
  variables, so the weak formulation of the Cauchy problem is needed
  to give meaning to the boundary condition in \eqref{eq:main}.

\begin{thm}[Well-posedness]\label{thm:existence}
  Let  $r\in L^{\infty}_{loc}(\Er^N\times\overline{\Theta})$ 
  be bounded from above (not necessarily periodic in $x$),
  and let $u_0\in C^0(\Er^N\times\overline\Theta)\cap L^\infty(\Er^N;L^1(\Theta))$ be nonnegative and bounded.
  Then, there exists a unique 
  solution $u$ to the Cauchy problem~\eqref{eq:main}-\eqref{eq:IC}.
  Moreover, $u$ and $\rho$ are bounded and nonnegative, and if $u_0\not\equiv 0$, then they are strictly positive.
\end{thm}

\subsection{Long-time behaviour of the solution}\label{ss:intro_ltb}

The long-time behaviour of solutions
of the Fisher-KPP equation is characterised by the sign of the principal
eigenvalue of the linearised operator around the null state.
This is a well-known fact when the problem is set on a  smooth bounded domain,
see e.g.\ \cite{CanCos89}. This is also known if the operator is periodic,
thanks to \cite{BHR05}. In both cases, the principal
eigenvalue is provided by the classical Krein-Rutman theory.
The characterisation has been extended in \cite{BHR07} to the non-compact setting 
by using the notion of the {\em generalised principal eigenvalue} 
inspired by~\cite{BNV94}.
It has then been employed in several different frameworks, see \eg \cite{Gir23,GirMaz24,HLR21,BCV16,ABR17,ACR13}.
This notion is required to study our model when
$\Theta=\Er^P$. 

More precisely, we call $\mcl[r]$
the linearised operator associated with~\eqref{eq:main}, that is
\[\mcl[r]\phi:=\Delta_x\phi+\Delta_\theta\phi+r(x,\theta)\phi.\]
We emphasise that $\mcl[r]$ is a local operator.

\begin{dfi}\label{dfi:eigenvalue}
  The generalised principal eigenvalue of $\mcl[r]$ is
  \[\lambda[r]:=\inf\acc{\lambda\in\Er\tq 
    \exists\phi>0,\ \mcl[r]\phi\loq\lambda\phi \text{ in }\Er^N\times\Theta,\ \nu\cdot\nabla_{\theta}\phi\goq 0\text{ on $\Er^N\times\dr\Theta$}}.\]
\end{dfi}
In the above definition, and throughout the paper, the functions $\phi$
are understood to belong to $W^2_{p,loc}(\Er^N\times\overline\Theta)$ 
for some
$p>1+N+P$, hence to be $C^1$ up to the boundary
by Morrey's embedding, so the boundary condition holds in the classical sense.

Definition \ref{dfi:eigenvalue} 
applies to both cases $\Theta$ bounded and $\Theta=\Er^P$
(in the latter, the boundary condition is ignored, as usual).
We point out that, even in the case of bounded~$\Theta$, the domain of
the functions $\phi$ on which the operator 
acts is not compact, since we do not assume $\phi$ to be periodic in $x$.
However, despite the lack of compactness, it turns out that $\lambda[r]$
is indeed a  {\em generalised principal eigenvalue}
for the operator $\mcl[r]$ acting
 on spatially-periodic functions, in the sense that it admits a positive 
 eigenfunction in that space.

 \begin{ppn}\label{ppn:pev}
 Let $r$ satisfy~\eqref{eq:general_assumptions} and let $\lambda[r]$
 be the generalised principal eigenvalue of $\mcl[r]$
 given by Definition \ref{dfi:eigenvalue}.
 There exists a solution $\varphi$ of the problem
\begin{equation}\label{eq:pb_vp}
    \left\{
  \begin{aligned}
    &\mcl[r]\varphi(x,\theta)=\lambda[r]\varphi(x,\theta), &(x,\theta)\in\Er^N\times\Theta,\\
    &\nu\cdot\nabla_{\theta}\varphi(x,\theta)=0, &\text{ $(x,\theta)\in\Er^N\times\dr\Theta$},\\
    &\varphi(x,\theta)>0, &(x,\theta)\in\Er^N\times\Theta,\\
    &\varphi\text{ is $1$-periodic in $x$}.
  \end{aligned}
  \right.
  \end{equation}
  \end{ppn}
  
We call {\em any} solution $\varphi$ of~\eqref{eq:pb_vp} a {\em principal eigenfunction} associated with $\lambda[r]$. 
When $\Theta$ is bounded, the Krein-Rutman theory \cite{KreRut48}
provides the existence, uniqueness and simplicity of
the principal eigenvalue for problem~\eqref{eq:pb_vp},
which furthermore can be expressed through the Rayleigh quotient;
in such a case, Proposition \ref{ppn:pev} asserts that this principal eigenvalue 
is characterised
by the formula in Definition \ref{dfi:eigenvalue}.
In the case where~$\Theta$ is unbounded, the Krein-Rutman theory does not apply
and indeed neither simplicity nor uniqueness hold for the eigenproblem~\eqref{eq:pb_vp}, and there is actually a half-line of eigenvalues
admitting positive eigenfunctions, see \cite{BerRos09}.
This is why we call $\phi$ in~\eqref{eq:pb_vp} ``a'' --~rather than ``the''~-- 
principal eigenfunction.

Making use of some results concerning the generalised principal eigenvalue,
we are able to establish the ``standard'' extinction/persistence dichotomy for the
the long-time behaviour of the solution for our model.

\begin{thm}\label{thm:ltb} 
  Let $r$ satisfy~\eqref{eq:general_assumptions}-\eqref{eq:r<0} 
  and let $\lambda[r]$ be the generalised principal eigenvalue of $\mcl[r]$.
  Let $u$ be the 
  solution to~\eqref{eq:main}-\eqref{eq:IC} with 
  $u_0\in C^0(\Er^N\times\overline\Theta)\cap L^\infty(\Er^N;L^1(\Theta))$ bounded, nonnegative and  not identically equal to $0$.
  \begin{itemize}
  \item If $\lambda[r]<0$, and in addition $u_0$ is compactly supported,
  then the population gets extinct, in the sense that
    \[\lim_{t\to+\infty}u(t,x,\theta)=0\quad\text{ uniformly in $(x,\theta)\in\Er^N\times\Theta$}.\]
  \item If $\lambda[r]>0$ then the population persists, in the sense that
    \begin{equation*}
    \liminf_{t\to+\infty}\cro{\sup_{x\in\Er^N}\rho(t,x)}>0.
    \end{equation*}
  \end{itemize}
\end{thm}

  We point out that the persistence property in the case $\lambda[r]>0$
  does not provide us with   a complete description of the long-time behaviour of the solution, but it just ensures that it stays bounded from below away from zero.
  In particular, we do not know whether the solution converges towards a positive stationary state (as in the case without phenotype~\eqref{eq:kpp_hetero} or in the spatially homogeneous case with phenotype~\eqref{eq:main_homo_intro}); 
  for instance, we cannot rule out the 
  possibility that the solution approaches
  a time-periodic cycle, which is a behaviour observed in some nonlocal
  reaction-diffusion equations, see e.g.~\cite{Toy}.
  
  As a matter of fact, the hypothesis that $u_0$ is compactly supported 
  in the persistence result is not needed in the case $\Theta$ bounded,
  cf.~Remark~\ref{rem:u0} below.

   Finally, the borderline case $\lambda[r]=0$ is not treated in this work. In some related reaction-diffusion models, it is known that 
   when the principal eigenvalue of the linearised operator is zero, 
   the population gets extinct. 
   This is the case when the phenotype is not taken into account \cite{BHR05}, or when the phenotype is taken into account through a quadratically decreasing Fisher geometric model 
   and the spatial structure is discrete, with the spatial variable limited to two possible values \cite{HLR21}.

\subsection{Influence of the heterogeneity of the environment and of the diffusivity}\label{ss:intro_period}

 We introduce two positive parameters $L,d$ in our model
by replacing the first equation in~\eqref{eq:main} by either
\begin{equation*}
\dr_tu=\Delta_xu+\Delta_\theta u+u\pth{r(x/L,\theta)-\rho(t,x)},
\qquad t> 0,\  (x,\theta)\in\Er^N\times\Theta,
\end{equation*}
or
\begin{equation*}
\dr_tu=d\Delta_xu+\Delta_\theta u+u\pth{r(x,\theta)-\rho(t,x)},
\qquad  t> 0,\  (x,\theta)\in\Er^N\times\Theta.
\end{equation*}
The quantity $1/L$ represents the frequency of oscillations of the landscape,  
while the diffusion coefficient $d$ accounts for the mobility of individuals.
Observe that one can pass from one of the above equations to the other
by a simple change of variable.
It is clear that all our previously stated
results hold true if one replaces the $1$-periodicity in $x$ with the $L$-periodicity,
hence the same is true for the model with diffusivity $d$.

Let us write for short 
\[r^L(x,\theta):=r\pth{\frac{x}{L},\theta}.\]
Thus, extinction or persistence for the model with $L$-periodicity
are determined by the sign of
the generalised principal eigenvalue $\lambda[r^L]$ given by
Definition \ref{dfi:eigenvalue}, while 
for the model with spatial diffusivity $d$ they are determined 
by the sign of the generalised principal eigenvalue $\lambda_d[r]$
of the linearised operator 
\[\mcl_d[r]:=d\Delta_x+\Delta_{\theta}+r.\]
We derive the following.

\begin{ppn}\label{ppn:deriv_lambda}
  Let $r$ satisfy the assumptions~\eqref{eq:general_assumptions}-\eqref{eq:r<0}
  and let $\lambda[r^L]$ and $\lambda_d[r]$ denote the generalised principal eigenvalue 
  of the linearised operator $\mcl[r^L]$ and $\mcl_d[r]$ respectively.
  \begin{enumerate}[label=$(\roman*)$]
  \item The function $L\mapsto\lambda[r^L]$ is continuous and nondecreasing on $\ioo{0,+\infty}$. 
  \item The function $d\mapsto\lambda_d[r]$ is
  continuous and nonincreasing on $\ioo{0,+\infty}$.
  \end{enumerate}
%
\end{ppn}

This result, combined with Theorem \ref{thm:ltb}, shows that
persistence becomes harder for the population if either
the frequency of oscillations of the landscape or the mobility of individuals 
increase.
This means that habitat fragmentation is armful for the species, and the same is true for 
the mobility. While the first fact is rather intuitive, the second one can be 
surprising, because one 
might expect that a faster diffusivity could provide higher chances of
adaptation. Our result shows that in our context, the opposite is true.
From a modelling perspective, 
this is in accordance with the conclusion of
\cite[Proposition {2.7.$i)$}]{HLR21}, where a spatially discrete system is considered
and an analogous monotonicity property with respect to the migration rate is obtained.

\subsection{A variant with heterogeneous diffusion}

For simplicity, we have 
considered that only the fitness of an individual is affected by its position and its phenotype. We may also deem it relevant to assume that the mobility and the mutation rate are anisotropic and depend on the position and the phenotype of the individual. This 
leads to the following extension of \eqref{eq:main}:
\begin{equation*}
  \begin{cases}
    \dr_tu=\nabla_x\cdot(a(x,\theta)\nabla_xu)+\nabla_\theta \cdot(\mu(x,\theta)
    \nabla_{\theta}u)\\[2pt]
    \ \esp+u\pth{r(x,\theta)-
      \displaystyle\int_\Theta{u}(t,x,\sigma)\de \sigma},
      &t>0,\ (x,\theta)\in \Er^N\times\Theta,\\[5pt]
    \nu\cdot(\mu(x,\theta)\nabla_{\theta}u) =0, &t>0,\ (x,\theta)\in\Er^N\times\dr\Theta,\\[5pt]
    u(0,x,\theta)=u_0(x,\theta), &(x,\theta)\in \Er^N\times\Theta,
  \end{cases}
\end{equation*}
where the functions $a,\mu$ are ranged in two sets
of uniformly positive definite, symmetric matrices of size $N\times N$ and $P\times P$
respectively,
and are $1$-periodic in the $x$ variable.
One may check that our arguments remain valid in this setting: the problem is well-posed, the long-time behaviour is described by the sign of the generalised principal eigenvalue 
(cf.~Theorems~\ref{thm:existence} and \ref{thm:ltb}), and the latter 
is continuous and nonincreasing with respect to the frequency of oscillations of the habitat and 
to the amplitude of the spatial diffusivity (cf. Proposition~\ref{ppn:deriv_lambda}). 

\section{Well-posedness of the model (Theorem~\ref{thm:existence})}\label{s:existence}

\begin{proof}[Proof of Theorem~\ref{thm:existence}]
  
We~derive the well-posedness of the Cauchy problem 
on the parabolic domain $[0,T]\times \Er^N\times\Theta$, with arbitrary $T>0$. By uniqueness, this will 
  provide us with a solution on the entire 
  $[0,+\infty)\times \Er^N\times\Theta$.
The proof is based on a contraction argument. To make it work, we start by deriving a priori bounds on $u$ and $\rho$.

In the sequel, the functions $r$ and $u_0$
fulfil the hypotheses of the theorem.

 \paragraph{Step 1. Global boundedness of $u$ and $\rho$.}
 
Let $u:[0,+\infty)\times\Er^N\times\Theta\to\Er$ be a solution to 
\eqref{eq:main}-\eqref{eq:IC} (in the weak sense, i.e., 
$u\in L^{\infty}((0,T)\times\Er^N\times\Theta)\cap L^{\infty} ((0,T)\times\Er^N;L^1(\Theta))$ for all $T>0$
and~\eqref{eq:weak_formulation} holds). 
    Since $u,\rho$ are bounded, locally in time, we can apply the parabolic maximum principle to \eqref{eq:main}, treating $\rho$ as a datum, and infer that $u(t,x,\theta)\geq0$ for all $t>0$, $(x,\theta)\in\Er^N\times\Theta$.
    Moreover, the strict inequality holds as soon as $u_0\not\equiv0$, owing to the 
    strong maximum principle.
  Furthermore, it holds that 
  \begin{equation}\label{eq:sub-heat}
  \dr_tu\loq\Delta_xu+\Delta_\theta u+
  (\bar r-\rho)\, u,\qquad t> 0,\  (x,\theta)\in\Er^N\times\Theta,
  \end{equation}
  where $\bar r:=\sup_{\Er^N\times\Theta}r$, which is finite by assumption.
  Again, this is understood in the weak sense, namely, the inequality 
  ``$\loq$'' holds in~\eqref{eq:weak_formulation} for nonnegative $\phi$'s. 
  Restricting to $\phi$ independent of $\theta$ in this 
  weak formulation,
  and taking first the integral on $\theta$ over $\Theta$ by Fubini's theorem,
  one finds that the function $\rho$ satisfies
  in the weak sense
  \begin{equation*}
      \dr_t\rho\leq \Delta_x\rho+(\bar r-\rho)\rho,
      \qquad t> 0,\ x\in\Er^N,
     \end{equation*}
together with the initial condition 
  $\rho(0,x)=\int_{\Theta} u_0(x,\sigma)\de \sigma$.
  From this, one gets, again by comparison,
  \begin{equation}\label{eq:estimate_rho}
    \forall t>0,\quad\|\rho(t,\cdot)\|_\infty\leq A:= \max\{\|\rho(0,\cdot)\|_\infty,\bar r \}.
  \end{equation}

  Last, we show that $u$ is globally bounded.
  For this, we treat the term $\rho$ as a given coefficient, that we know being bounded, of the equation in \eqref{eq:main}, which is then seen
  as a local linear parabolic equation.
  
  We first deal with the case $\Theta$ bounded. We use the parabolic Harnack inequality.
  In order to rule out the boundary of $\Theta$, we
  use the ${C}^2$-regularity of $\dr\Theta$ to
  extend the solution $u$ in the $\theta$ variable
  by orthogonal reflection with respect to 
  $\partial\Theta$ and get a function defined for all $\theta$ on an open
  set $\widetilde\Theta$ containing~$\overline\Theta$.
  The extended function is a solution of a parabolic equation
  of the same type as the one in~\eqref{eq:main}, with
  $\Delta_\theta$ extended outside $\Theta$ by an elliptic operator with coefficients
  depending on~$\Theta$, hence regular (Lipschitz-continuous).
  See \eg the proof of Theorem~3.1 in~\cite{BerRos09}.
  This allows us to apply the parabolic interior  Harnack inequality, 
  cf. \cite[Theorem~1.1]{KS80}, which,
  for any $R>0$,
  gives a constant $C_R>0$ such that
  \begin{equation}\label{Harnack-bdd}
  \forall t\geq 1,\ \forall x\in\Er^N,\qquad
  \sup_{B_{\Er^N}(x,R)\times\Theta}u(t,\cdot,\cdot)\leq 
  C_R\inf_{B_{\Er^N}(x,R)\times\Theta}u(t+1,\cdot,\cdot),
  \end{equation}
  where $B_{\Er^N}(x,R)$ is the open ball in $\Er^N$ with centre $x$ and radius $R$.
  We emphasise that the constant $C_R$ is
  independent of $t$ and $x$, because the coefficients of the equation (including~$\rho$)
  are uniformly bounded. Since 
  $\rho(t+1,x)\geq|\Theta|\inf_{\Theta}u(t+1,x,\cdot)$,
   with $|\Theta|$ denoting the measure of $\Theta$, we deduce in particular
   $$\forall t\geq 1,\ \forall x\in\Er^N,\qquad
  \sup_{\Theta}u(t,x,\cdot)\leq 
  C_1\inf_{\Theta}u(t+1,x,\cdot)\leq \frac{C_1}{|\Theta|}\rho(t+1,x),$$
  hence the global boundedness of $u$ follows from the one of $\rho$.

    Now, let us deal with the case $\Theta=\Er^P$.
    We see $u$ as a subsolution to the equation~\eqref{eq:sub-heat} and we apply the local maximum principle, see e.g.~\cite[Theorem~7.36]{Lie96}:
    there exists a constant $C>0$ such that for any $t\geq 2$, $x_0\in\Er^N$ and $\theta_0\in\Er^P$, there holds:
    \[\sup_{(t-1/2,\, t)\times B_{\Er^N}(x_0,1/2)\times B_{\Er^P}(\theta_0,1/2)}u\leq 
    C\int_{t-1}^t\int_{B_{\Er^N}(x_0,1)}\int_{B_{\Er^P}(\theta_0,1)}u.\]
  We emphasise that the constant $C$ is
  independent of $t$ and $x_0$, because the coefficients of the equation~\eqref{eq:sub-heat} (including~$\rho$)
  are uniformly bounded.
  We obtain in particular:
  \[u(t,x_0,\theta_0)\leq C\int_{t-1}^t\int_{B_{\Er^N}(x_0,1)}\rho\leq C|B_{\Er^N}(x_0,1)|\,\sup\rho,\]
  where $|B_{\Er^N}(x_0,1)|$ denotes the volume of $B_{\Er^N}(x_0,1)$.
  Then, $\rho$ being uniformly bounded, we conclude that $u$ is uniformly bounded
  too. 

  \paragraph{Step 2. Contraction argument.} 
  
Now, we fix an arbitrary $T>0$ and we take $\tau\in(0,T)$ to be determined later.
For a given function $w:[0,\tau]\times \Er^N\to L^1(\Theta)$,
we set 
\[\rho[w](t,x):=\int_{\Theta}w(t,x,\sigma)\de\sigma.\]
Then, we consider the set of functions $K$ defined by
\[K:=\acc{
w\in
L^\infty([0,\tau]\times\Er^N;L^1(\Theta))\tq 
0\loq \rho[w](t,x)\loq 2A e^{\bar r\, T}
\ \text{ for }(t,x)\in [0,\tau]\times\Er^N} .\]
Next, for a given initial datum $\varpi\in L^\infty(\Er^N\times\Theta)\cap
L^\infty(\Er^N;L^1(\Theta))$ 
satisfying 
\begin{equation}\label{eq:rho0}
\varpi\goq 0 \quad\text{in }\Er^N\times\Theta,\qquad
\text{ess}\sup_{x\in\Er^N}\|\varpi(x,\cdot)\|_{L^1(\Theta)}\loq A e^{\bar r\, T}\,,
  \end{equation}
and for given $w\in K$, we consider the following problem:
\[
(E_{w})\esp
  \left\{
  \begin{aligned}
    \dr_tu&=\Delta_xu+\Delta_\theta u+(r(x,\theta)-\rho[w](t,x))u,
    \hfill
    &(t,x,\theta)\in (0,\tau)\times\Er^N\times\Theta,\\
    \nu\cdot\nabla_{\theta} u&=0, &(t,x,\theta)\in (0,\tau)\times\Er^N\times\dr\Theta,\\
    u(0,x,\theta)&=\varpi(x,\theta), &(x,\theta)\in\Er^N\times\Theta.
  \end{aligned}
  \right.
  \]
  Since $r$ is bounded from above and $\rho[w]$ is bounded, 
  the comparison principle holds for the above problem, and 
  it is standard that it 
  admits a unique bounded weak solution (in the sense precised at the beginning
  of Section \ref{sec:well}),
  which we call $u^w$. 
  It also follows that $u^w\in W^{1,2}_{p,loc}((0,\tau]\times\Er^N\times\overline{\Theta})$ and that it is nonnegative,
  and thus by the same argument as in the Step 1, the function
  $\rho[u^w](t,x)$ fulfils in the weak sense
   \begin{equation*}
      \dr_t\rho[u^w]\leq \Delta_x\rho[u^w]+(\bar r-\rho[w])\rho[u^w],
      \qquad t\in(0,\tau),\ x\in\Er^N.
     \end{equation*}
 Then, owing to  \eqref{eq:rho0}, we derive
  \begin{equation*}
  \forall t\in[0,\tau],\,\forall x\in\Er^N,\qquad \rho[u^w](t,x)\loq 
  Ae^{\bar r\, (T+t)}\,.
  \end{equation*}
  It follows that~$u^w$ belongs to $K$ provided $e^{\bar r \, \tau}\loq 2 $.

  Next, 
   take $w_1,w_2\in K$. Let us write for short
   $u_1:=u^{w_1}$, $u_2:=u^{w_2}$ and $\rho_1:=\rho[w_1]$, $\rho_2:=\rho[w_2]$. 
   The function $v:=u_1-u_2$ vanishes identically at $t=0$ and satisfies
   \begin{align}
  \dr_tv -\Delta_x v-\Delta_{\theta}v&=rv-\rho_1v-(\rho_1-\rho_2)u_2\nonumber\\
  & \loq (r-\rho_1)v+
  \|w_1-w_2\|_{L^\infty([0,\tau]\times\Er^N;L^1(\Theta))}
  u_2\,.\label{eq:v_subsolution}
   \end{align}
    Let $\bar v$ be the solution of the equation
    \[\dr_t\bar v -\Delta_x \bar v-\Delta_{\theta}\bar v= \bar r\,\bar v
    +\|w_1-w_2\|_{L^\infty([0,\tau]\times\Er^N;L^1(\Theta))}
    u_2,
    \qquad (t,x,\theta)\in (0,\tau)\times\Er^N\times\Theta,\]
  with Neumann boundary conditions and initial datum identically equal to $0$.
    On the one hand,~$\bar v$ is nonnegative, so
    we have $\bar r\bar v\goq (r-\rho_1)\bar v$; thus $\bar v$ is a supersolution of
    the equation~\eqref{eq:v_subsolution}, for which $v$ is a subsolution. 
    We deduce by comparison
    that $v\loq\bar v$. On the other hand, 
    in the case $\Theta$ bounded, as done before, 
    we can integrate the equation for $\bar v$ on $\Theta$. 
    Hence, since
    $u_2\in K$, we get 
    \[\begin{split}
    \dr_t\rho[\bar v] -\Delta_x \rho[\bar v] &= \bar r\,\rho[\bar v]
    +\|w_1-w_2\|_{L^\infty([0,\tau]\times\Er^N;L^1(\Theta))}
    \rho[u_2]\\
    &\loq \bar r\,\rho[\bar v]
    +\|w_1-w_2\|_{L^\infty([0,\tau]\times\Er^N;L^1(\Theta))}
    \times 2A e^{\bar r\, T}\,.
    \end{split}\]
    It follows that, for $\tau$ sufficiently small, only depending on $\bar r$,
    $A$ and $T$, it holds  
    $$ \forall t\in[0,\tau],\,\forall x\in\Er^N,\qquad 
    \rho[\bar v](t,x)\loq \frac12\|w_1-w_2\|_{L^\infty([0,\tau]\times\Er^N;L^1(\Theta))}\,.$$
    Instead, in the case $\Theta=\Er^P$,
    one can explicitly compute the solution $\bar v$ through the heat kernel 
    and check 
    that the above estimate holds true, again for 
    $\tau$ depending on $\bar r$, $A$ and $T$.
    This gives a one-side inequality for $\rho[ v]$. The same argument 
    for $-\rho[v]=\rho[-v]$ follows reverting the roles of $u_1$ and $u_2$.
    In conclusion, we have shown that the mapping $w\mapsto u^w$
    is a contraction from the closed set $K$ into itself, with respect to
    the $L^\infty([0,\tau]\times\Er^N;L^1(\Theta))$
    norm, for any $\tau\in(0,\tau_0]$, where $\tau_0>0$ is sufficiently small and depends only on $\bar r$, $A$ and $T$.
    Therefore, $w\mapsto u^w$ has a unique fixed point, which implies that~\eqref{eq:main} has a
    unique solution on $(0,\tau)\times\Er^N\times\Theta$ for any 
    $\tau\in(0,\tau_0]$ and any initial condition $\varpi$ satisfying~\eqref{eq:rho0}.

\paragraph{Step 3. Conclusion.}
    Step 2 provides us with a unique solution $u$ to
    problem \eqref{eq:main} $(0,\tau_0)\times\Er^N\times\Theta$
    with initial datum $u_0$. 
    Take $\eps\in\ioo{0,1/2}$ and consider 
    the new initial datum $\varpi=u((1-\eps)\tau_0,\cdot,\cdot)$.
    By~\eqref{eq:estimate_rho}, $\varpi$ satisfies~\eqref{eq:rho0},
    hence by Step 2 the corresponding problem $E_w$ 
    admits a unique solution on $(0,\tau_0)\times\Er^N\times\Theta$, 
    which overlaps with $u(\cdot+(1-\eps)\tau_0,\cdot,\cdot)$ on 
    $[0,\eps\tau_0)\times\Er^N\times\Theta$. 
    We may thus extend $u$ to a solution of \eqref{eq:main} for $t\in\ifo{0,(2-\eps)\tau_0}$. 
    Iterating this argument, we have that $u$ solves \eqref{eq:main} for all times
    $t\in(0,T]$, hence for all $t>0$, by the arbitrariness of $T$.
    We further know by Step 1 that $u$ is globally bounded, together with its
    associated function~$\rho$.
\end{proof}

\section{Properties of the principal eigenfunction and eigenvalue}

Before studying the long-time behaviour of the solution, we 
derive 
some properties of the generalised principal eigenvalue $\lambda[r]$ 
of Definition \ref{dfi:eigenvalue}.

We start with the case $\Theta$ bounded.
We derive a result containing two statements: 
the first one follows from the arguments of \cite{BerRos09} and provides us
with an approximation of 
the principal eigenvalue~$\lambda[r]$ by {principal eigenvalues} 
in the domains $B_{\Er^N}(0,R)\times\Theta$ under mixed Dirichlet-Neumann
boundary conditions.
The second statement asserts that $\lambda[r]$ 
coincides with the standard periodic principal eigenvalue on $\Er^N\times\Theta$.
This implies in particular that Proposition \ref{ppn:pev} holds in the case
$\Theta$ bounded.

\begin{ppn}\label{ppn:br_bounded}
  Assume $\Theta$ is bounded and has ${C}^2$ boundary. Let $r$ satisfy~\eqref{eq:general_assumptions}.
  \begin{enumerate}[label=$(\roman*)$]
  \item   For all $R>0$, 
   the mixed eigenproblem
  \begin{equation*}
  \begin{cases}
    \mcl[r]\varphi_R^{m}=\lambda_R^{m}\varphi_R^{m} & \text{in }B_{\Er^N}(0,R)\times\Theta,\\
    \nu\cdot\nabla_{\theta}\varphi_R^{m}=0 & \text{over $B_{\Er^N}(0,R)\times\dr \Theta$},\\
    \varphi_R^{m}=0 & \text{over $(\partial B_{\Er^N}(0,R))\times \Theta$},\\
    \varphi_R^{m}>0 & \text{in }B_{\Er^N}(0,R)\times\Theta,
  \end{cases}
  \end{equation*}
  admits a unique solution $\lambda^m_R\in\Er$ and a unique (up to a scalar multiple) generalised solution $\varphi_R^{m}\in W^{1}_2(B_{\Er^N}(0,R)\times\Theta)\cap 
  C^0( \overline{B_{\Er^N}(0,R)\times\Theta})$.
  Moreover, the function $R\mapsto\lambda^{m}_R$ is 
  increasing and satisfies $\lambda^{m}_R\to \lambda[r]$
  as $R\to+\infty$.
Finally, the following Rayleigh formula~holds:
  \[\lambda_R^m=\max_{\phi\in X_R}\pth{\frac{-\displaystyle\int_{B_{\Er^N}(0,R)\times\Theta}\abs{\nabla \phi}^2+\int_{B_{\Er^N}(0,R)\times\Theta}
      r\phi^2}{\displaystyle\int_{B_{\Er^N}(0,R)\times\Theta}\phi^2}},
  \]
  where
  \[X_R:=\acc{\phi\in W^{1}_2(B_{\Er^N}(0,R)\times\Theta)\cap C^0(  \overline{B_{\Er^N}(0,R)\times\Theta})\ /\ \phi\equiv 0 \text{ over } \dr B_{\Er^N}(0,R)\times\Theta}.\]
  
\item The problem 
  \begin{equation*}
    \left\{
    \begin{aligned}
      &\mcl[r]\varphi(x,\theta)=\lambda\varphi(x,\theta) &\text{in $\Er^N\times\Theta$},\\
      &\nu\cdot\nabla_{\theta}\varphi(x,\theta)=0 &\text{ over $\Er^N\times\dr\Theta$},\\
      &\varphi(x,\theta)>0 &\text{in $\Er^N\times\Theta$},\\
      &\varphi\text{ is $1$-periodic in $x$},
    \end{aligned}
    \right.
  \end{equation*}
  admits as a solution the unique eigenvalue $\lambda=\lambda[r]$ and 
  a unique (up to a scalar multiple) eigenfunction 
  $\varphi\in W^{1}_{2,loc}(\Er^N\times{\overline\Theta})$.
  In particular, the generalised principal eigenvalue given by Definition \ref{dfi:eigenvalue} 
  coincides with the classical $x$-periodic principal eigenvalue of the above problem.

    \end{enumerate}
\end{ppn}

\begin{proof}
  \textbf{Statement $(i)$.}  We point out that the eigenproblem stated in $(i)$ is not
  completely standard due to the presence of the ``corners'' 
  $\partial B_{\Er^N}(0,R))\times\partial \Theta$.
  However, the existence and uniqueness up to a scalar multiple of $\varphi_R^{m}$ is derived in 
  \cite[Theorem~3.1]{BerRos09}  in dimension $N=1$. The arguments, which rely on the solvability, in the generalised sense,
  of the mixed boundary value problem
  (see e.g.\ the Notes to Chapter 8 in \cite{GilTru15}) and on the construction of some barriers to get continuity up to
  the corners, holds true without modification in the higher dimensional case.
  The monotonicity and the 
  convergence of $\lambda_R^{m}$ to $\lambda[r]$ are
  due to \cite[Proposition 1]{BerRos09} (again in dimension $N=1$,
  but the argument works in arbitrary dimension). 

  \medskip

  We now show the Rayleigh formula for $\lambda^m_R$.
    Let $W_{2,0}^{1}(B_{\Er^N}(0,R)\times\overline{\Theta})$ denote the closure in $W^{1}_2(B_{\Er^N}(0,R)\times\Theta)$ of
    \[\acc{\phi\in{C}^1(B_{\Er^N}(0,R)\times\overline{\Theta})\ / \ \phi\equiv 0 \text{ over } \dr B_{\Er^N}(0,R)\times\Theta}.\]
    By the classical elliptic theory (see \eg the Notes to Chapter 8 in~\cite{GilTru15}), 
    for all $f\in L^2(B_{\Er^N}(0,R)\times\Theta)$, there exists a unique weak solution $\phi\in W_{2,0}^{1}(B_{\Er^N}(0,R)\times\overline{\Theta})$ to the mixed boundary value problem
  \[
  \left\{
  \begin{aligned}
    &(\Delta_x+\Delta_{\theta})\phi+(r-\bar r)\phi=f &\text{in $B_{\Er^N}(0,R)\times\Theta$},\\
    &\nu\cdot\nabla\phi=0&\text{over $B_{\Er^N}(0,R)\times\dr\Theta$},\\
    &\phi=0&\text{over $\dr B_{\Er^N}(0,R)\times\Theta$},
  \end{aligned}
  \right.
  \]
  where $\bar r=\sup r$.
   Let $T\,:\;L^2(B_{\Er^N}(0,R)\times\Theta)\to L^2(B_{\Er^N}(0,R)\times\Theta)$ be the
    operator which assigns to $f$ 
    the solution $\phi$.
    Owing to the ${C}^2$ regularity of $\Theta$ 
    and to the Neumann boundary condition,
    one may extend $\phi$ by reflection 
    to a larger domain $B_{\Er^N}(0,R)\times\widetilde\Theta$ with 
    $\overline\Theta\subset\widetilde\Theta$, which satisfies there a uniformly elliptic equation. 
    This allows one to apply
    standard elliptic estimates and a compact injection theorem to
    infer that the operator $T$ is compact.
    See \eg the argument of Step~2 of the proof of~\cite[Theorem~3.1]{BerRos09}.
  Thus, the operator $T$ being symmetric, 
  there exists an orthonormal basis composed of eigenvectors of $T$.
  Those eigenvectors are also eigenvectors of the mixed eigenproblem 
  \begin{equation*}
  \begin{cases}
    \mcl[r]\phi=\lambda_k\phi & \text{in }B_{\Er^N}(0,R)\times\Theta,\\
    \nu\cdot\nabla_{\theta}\phi=0 & \text{over $B_{\Er^N}(0,R)\times\dr \Theta$},\\
    \phi=0 & \text{over $(\partial B_{\Er^N}(0,R))\times \Theta$},
    \end{cases}
  \end{equation*}
  for some eigenvalues $\lambda_k\in\Er$.
  Therefore, the Rayleigh formula holds for the mixed eigenproblem, see \eg the proof of  
  \cite[Theorem~6.5.2]{Eva98}:
  \[\lambda_R^m=\max_{\phi\in W_{2,0}^{1}(B_{\Er^N}(0,R)\times\overline{\Theta})}\pth{\frac{-\displaystyle\int_{\Er^N\times\Theta}\abs{\nabla \phi}^2+\int_{\Er^N\times\Theta}
      r\phi^2}{\displaystyle\int_{\Er^N\times\Theta}\phi^2}}.\] 
  We point out that the $\max$ is reached by the principal eigenfunction $\varphi_R^m$, which belongs to the space~$X_R$ defined in
  the statement of the proposition. 
  Therefore, we may replace $W_{2,0}^{1}(B_{\Er^N}(0,R)\times\overline{\Theta})$ by 
  the subspace $X_R$ in the variational formula.

  \medskip
  
  \noindent\textbf{Statement $(ii)$.} The classical Krein-Rutman theory~\cite{KreRut48} (see \eg the proof of~\cite[Theorem~5.1]{BerRos09}) yields the
  uniqueness and simplicity of the Neumann/periodic eigenproblem of statement~$(ii)$.
  Namely, the problem admits a unique eigenvalue $\lambda$ and 
  a unique (up to a scalar multiple) eigenfunction $\varphi$.
  What we need to show is that $\lambda=\lambda[r]$. 
  
The inequality $\lambda[r]\leq\lambda$ directly follows from 
Definition \ref{dfi:eigenvalue}. In order to show the reverse inequality, we make
use of the Rayleigh formula for $\lambda_R^m$.
Consider a family of cutoff functions $(\zeta_R)_{R>1}$
  in ${C}^{\infty}(\Er^N)$ satisfying $\zeta_R\equiv 1$ in $B_{\Er^N}(0,R-1)$, $\zeta_R\equiv 0$ 
  on $\partial B_{\Er^N}(0,R)$,
    $0\leq\zeta_R\leq 1$ and
  $\dabs{\nabla\zeta_R}_{\infty}\loq 2$.
  We then apply the Rayleigh formula  given in $(i)$ to the 
  test function $\phi=\varphi\zeta_R\in X_R$.
   Let us call for short $B_R:=B_{\Er^N}(0,R)$ and, for $R>1$,
      \[U_R:=\bigcup_{h\in\Zed^N\text{ s.t. } \mcc+\acc{h}\,\subset\, B_{R-1}}\mcc+\acc{h},\]
       where $\mcc:=[0,1]^N$ is the periodic cell.
  We get:
      \begin{align*}
        \lambda_R^m&\geq
        \frac{-\displaystyle\int_{U_R\times\Theta}\abs{\nabla \varphi}^2
        +\int_{U_R\times\Theta}r\varphi^2
        -\displaystyle\int_{(B_R\setminus U_R)\times\Theta}
        \big(\abs{\nabla (\varphi\zeta_R)}^2+ |r|\varphi^2\big)}{\displaystyle\int_{B_R\times\Theta}\varphi^2}\,.
       \end{align*}
       Using the $1$-periodicity of $\varphi$ and $r$, and noticing that
       $U_R$ contains exactly $|U_R|$ distinct cells $\mcc+\{h\}$, $h\in\Zed^N$,
       we obtain
       \begin{align*}
        \lambda_R^m&\geq
        \frac{|U_R|\pth{-\displaystyle\int_{\mcc\times\Theta}\abs{\nabla \varphi}^2
            +\int_{\mcc\times\Theta}r\varphi^2}-
            |B_R\setminus U_R|\times|\Theta|
            \big[\sup(|\nabla \varphi|+2\varphi)^2+|r|\varphi^2\big]}
               {\displaystyle|U_R|\int_{\mcc\times\Theta}\varphi^2
                 +\displaystyle|B_R\setminus U_R| \times|\Theta|\,\sup\varphi^2},
      \end{align*}    
      where $|\cdot|$ stands for the volume. 
      Therefore, since $\frac{|B_R\setminus U_R|}{|U_R|}\to 0$ as $R\to+\infty$, 
      letting $R\to+\infty$ in the above expression yields
      \[\lambda[r]=\lim_{R\to+\infty}\lambda_R^m\geq \frac{{-\displaystyle\int_{\mcc\times\Theta}\abs{\nabla \varphi}^2+\int_{\mcc\times\Theta}r\varphi^2}}
               {\displaystyle\int_{\mcc\times\Theta}\varphi^2}\,,\]
               and the right-hand side is equal to $\lambda$, thanks to the
               standard Rayleigh formula for $\lambda$ (or even by direct
               computation, integrating by parts).
    Thus, by $(i)$ we conclude $\lambda[r]\geq\lambda$. 
\end{proof}

Let us turn to the case $\Theta=\Er^P$.
Using again the results from \cite{BerRos09}, we now show
that the principal eigenvalue~$\lambda[r]$ of~$\mcl[r]$ 
coincides with the limit as $R\to+\infty$ of both the 
{Dirichlet principal eigenvalue} in the domain $B_{\Er^{N+P}}(0,R)$, and 
the $x$-periodic principal eigenvalue in the domain $\Er^N\times B_{\Er^P}(0,R)$.
These last two notions of principal eigenvalue are provided by the 
classical Krein-Rutman theory \cite{KreRut48}.
Namely, the eigenproblems  
\begin{equation}\label{eq:D_pevr}
  \begin{cases}
    \mcl[r]\varphi_R^{d}=\lambda_R^{d}\varphi_R^{d} & \text{in }B_{\Er^{N+P}}(0,R),\\
    \varphi_R^{d}=0 & \text{over $\partial B_{\Er^{N+P}}(0,R)$},\\
    \varphi_R^{d}>0 & \text{in $B_{\Er^{N+P}}(0,R)$},
  \end{cases}
  \end{equation}
and
  \begin{equation}\label{eq:periodic_pevr}
  \begin{cases}
    \mcl[r]\varphi_R=\lambda_R\varphi_R & \text{in }\Er^N\times B_{\Er^{P}}(0,R),\\
    \varphi_R=0 & \text{over $\Er^N\times\partial B_{\Er^{P}}(0,R)$},\\
    \varphi_R>0 & \text{in $\Er^N\times B_{\Er^{P}}(0,R)$},\\
    \varphi_R  \text{ is $1$-periodic in $x$},
  \end{cases}
  \end{equation}
admit a unique eigenpair $(\lambda^d_R,\varphi^d_R)$ and $(\lambda_R,\varphi_R)$ respectively (uniqueness for
$\varphi^d_R$ and $\varphi_R$ is understood up to a scalar multiple).

\begin{ppn}\label{ppn:cv_approximation_lambda}
  Assume $\Theta=\Er^P$ and let $r$ satisfy~\eqref{eq:general_assumptions}. Let $(\lambda^d_R,\varphi^d_R)$ and $(\lambda_R,\varphi_R)$
  be the eigenpairs of the problems \eqref{eq:D_pevr} and \eqref{eq:periodic_pevr}
  respectively. 
  \begin{enumerate}[label=$(\roman*)$]
  \item   The functions $R\mapsto\lambda^d_R$ and $R\mapsto\lambda_R$ are nondecreasing 
  and satisfy
  \[\lambda^d_R\to\lambda[r],\esp \lambda_R\to\lambda[r],
  \esp\text{as }\;R\to+\infty.\]
\item The principal eigenfunctions $\varphi_R$, normalised by 
$\varphi_R(0,0)=1$, converge
locally uniformly, as $R\to+\infty$, 
to a function $\varphi\in W^{2}_{p,loc}(\Er^N\times\Er^P)$ satisfying the periodic principal eigenproblem
\begin{equation*}
  \begin{cases}
    \mcl[r]{\varphi}={\lambda}{\varphi} & \text{in }\Er^N\times\Er^P,\\
    \varphi>0 & \text{in }\Er^N\times\Er^P,\\
    \varphi  \text{ is $1$-periodic in $x$.}    
  \end{cases}
\end{equation*}
  \end{enumerate}
\end{ppn}

\begin{proof}
  On the one hand, \cite[Proposition 4]{BerRos09} implies that $\lambda^d_R\to \lambda[r]$ as $R\to+\infty$ and that $R\mapsto\lambda^d_R$ is nondecreasing.
  On the other hand, \cite[Proposition 5]{BerRos09} and its proof imply that $R\mapsto\lambda_R$ is nondecreasing and that there exist ${\lambda}\in\Er$ and ${\varphi}\in W^{2}_{p,loc}(\Er^N\times\Er^P)$ such that $\lambda_R\to\lambda$, $\varphi_R\to\varphi$ locally uniformly, and
  \begin{equation*}
  \begin{cases}
    \mcl[r]{\varphi}={\lambda}{\varphi} & \text{in }\Er^N\times\Er^P,\\
    \varphi>0 & \text{in }\Er^N\times\Er^P,\\
    \varphi  \text{ is $1$-periodic in $x$.}    
  \end{cases}
  \end{equation*}
  By~\cite[Proposition~6]{BerRos09}, we also have $\lambda=\lambda[r]$~\footnote{\ We point out that although these results are stated in \cite{BerRos09} for bounded potential~$r$, the proofs rely on~\cite{BerRos15}, which only 
  requires $r$ to be locally bounded, 
  which is the case in our setting.}
  (in the notation
  of~\cite{BerRos09}, $\lambda=-\lambda_{1,l}$ and $\lambda[r]=-\lambda_1$).
  This concludes the proof.
\end{proof}

 Statements $(ii)$ of Propositions \ref{ppn:br_bounded} and \ref{ppn:cv_approximation_lambda}
yield Proposition \ref{ppn:pev} in both cases $\Theta$ bounded and unbounded.

\section{Long-time behaviour of the solution (Theorem \ref{thm:ltb})}
\label{s:ltb}


We start with a technical lemma which, together with hypothesis~\eqref{eq:r<0},
will allow us to translate the $L^1_{loc}(\Theta)$ convergence of $u$ to $0$ 
into $L^1(\Theta)$ convergence. This will be used to prove the extinction
result in the case $\Theta=\Er^P$ and,
curiously, also the persistence property.
In contrast with Harnack-type inequalities, the proof exploits the nonlinear 
term of the equation.

\begin{lem}\label{lem:rho<}
  Let $\Theta=\Er^P$ and let $r$ satisfy~\eqref{eq:general_assumptions}-\eqref{eq:r<0}. Let $u$ solve~\eqref{eq:main}-\eqref{eq:IC} with $u_0\in C^0(\Er^N\times\overline\Theta)\cap L^{\infty}(\Er^N;L^1(\Theta))$
  nonnegative and not identically equal to $0$.
    For given $t_0\goq 0$ and  $\tau,M>0$, call
    $$H:=\bar r\,\sup_{s\in(t_0,t_0+\tau),\; x\in\Er^N}
    \int_{\|\sigma\|\loq M}  u(s,x,\sigma)\de \sigma
  +\|\rho\|_\infty\,\sup_{x\in\Er^N,\;\|\theta\|> M}r^+(x,\theta),$$
  where $\bar r:=\displaystyle\sup_{\Er^N\times\Er^P} r$,
  and $r^+:=\max\{r,0\}$ denotes the positive part of $r$.
    Then it holds~that
    \begin{equation}\label{rho<H}
    \sup_{x\in\Er^N}\,\rho(t_0+\tau,x)\loq 
    \sqrt{H}\,\coth\left(\sqrt{H}\,\tau+
    \arcoth\pth{\frac{\|\rho\|_\infty+\sqrt{H}}{\sqrt{H}}}\right).
    \end{equation}
  Moreover, the expression at the right-hand side in~\eqref{rho<H} is increasing with respect to 
  $H$, for all~$\tau\goq0$. 
\end{lem}

\begin{proof}
    We recall that Theorem \ref{thm:existence} ensures that both $u$ and $\rho$
    are bounded.
    Take $t_0\goq 0$, $\tau,M>0$, and define $H$ as in the statement of the 
    lemma. 
    As in the proof of Theorem~\ref{thm:existence}, restricting to $\theta$-independent
    test functions in the weak formulation~\eqref{eq:weak_formulation}, and using Fubini's theorem, one finds that the function~$\rho$ satisfies in the weak sense
  \begin{equation}\label{eq:rho_pde}
  \dr_t\rho=\Delta_x\rho+\int_{\Er^P}r(x,\sigma)u(t,x,\sigma)\de \sigma-\rho^2,
    \qquad t>0,\ x\in\Er^N.
    \end{equation}
  For $t\in(t_0,t_0+\tau)$ and $x\in\Er^N$, we see that 
  \[\begin{split}
  \dr_t\rho-\Delta_x\rho &\loq 
     {\bar r}\int_{\|\sigma\|\loq M}  u(t,x,\sigma)\de \sigma
     +\|\rho\|_\infty\,\sup_{x\in\Er^N,\;\|\theta\|> M}r^+(x,\theta)-\rho^2
  \loq H-\rho^2.\\
  \end{split}\]
  Let $V$ be the solution of
  $V'=H-V^2$, with initial datum $V(t_0)=\|\rho\|_\infty+\sqrt{H}$. 
  We point out that $\coth'=1-\coth^2$; thus
  \[\forall t\geq 0,\quad
  V(t_0+t)=\sqrt{H}\coth\pth{\sqrt{H}t+\arcoth\pth{\frac{\|\rho\|_\infty+\sqrt{H}}{\sqrt{H}}}}.\]
  Moreover, $\rho(t_0,x)\loq V(t_0)$, so by the parabolic comparison principle, $\rho(t_0+t,x)\loq V(t_0+t)$ for 
  $t\in[0,\tau]$
  and $x\in\Er^N$. The inequality~\eqref{rho<H} is proved.

  Finally, Grönwall's inequality implies that 
  $V$ is strictly increasing with respect   to $H$. The proof is concluded.
\end{proof}

\begin{proof}[Proof of Theorem \ref{thm:ltb}, item 1]
  Assume that $\lambda[r]<0$.
  Let $\varphi>0$ be a generalised principal eigenfunction corresponding to $\lambda[r]$, i.e.~a solution of~\eqref{eq:pb_vp}.
  The function $e^{\lambda[r] t}\varphi(x,\theta)$
  is $1$-periodic in $x$ and it is a positive solution of the linear equation
  \[\dr_tv=\Delta_xv+\Delta_\theta v+r(x,\theta)v,\]
  and additionally fulfils the boundary condition $\nu\cdot\nabla_{\theta} v=0$
  if $\Theta$ is bounded. 
  The function $u$ is a subsolution of the above equation.
  Moreover, since $u_0$ has a compact support in $\Er^N\times\overline\Theta$,
and $\varphi>0$ on $\Er^N\times\overline\Theta$ (owing to the Hopf lemma in the case
  $\Theta$ bounded),
  there exists a large $a>0$ such that $a\varphi\goq u_0$.
  It follows from the comparison principle that $u(t,x,\theta)\loq a e^{\lambda[r] t}\varphi(x,\theta)$ for all $t\goq 0$, $(x,\theta)\in\Er^N\times\Theta$.
  This concludes the proof in the case $\Theta$ bounded.  When
  $\Theta=\Er^P$, we deduce that
 $u(t,x,\theta)\to 0$ as $t\to+\infty$, uniformly in 
 $x\in\Er^N$ and locally uniformly in $\theta\in\Er^P$.
 It then follows from Lemma \ref{lem:rho<} and assumption~\eqref{eq:r<0}
 that $\rho(t,x)\to 0$ as $t\to+\infty$, uniformly in 
 $x\in\Er^N$, which implies in turn, owing to the Harnack inequality,
 that $u(t,x,\theta)\to 0$ as $t\to+\infty$, uniformly in 
 $\theta\in\Er^P$ as well.
\end{proof}

\begin{req}\label{rem:u0}

The hypothesis that $u_0$ is compactly supported is only used in the 
above proof in order to ensure the inequality $a\varphi\goq u_0$
for large enough $a$. However, when $\Theta$ is bounded,
this property automatically holds because $\inf \varphi>0$
due to the periodicity in $x$ and the Hopf~lemma. 
Thus, when $\Theta$ is bounded, the assumption that $u_0$ has compact support
can be dropped.
When $\Theta=\Er^P$, again due to the periodicity in $x$ of $\varphi$, the assumption that $u_0$ has compact support can be weakened to the existence 
of some $A>0$ such that the support of $u_0$ is included in $\Er^N\times B_{\Er^P}(0,A)$.
\end{req}

We now turn to the proof of the persistence property. In this 
case we cannot neglect the nonlocal nonlinear term $-\rho u$, 
as we have done in the proof of extinction,
so we need to handle it in some way, despite it prevents the validity of the maximum principle. 
We treat separately the case where $\Theta$ is bounded and the case where $\Theta$ is unbounded.

\begin{proof}[Proof of Theorem \ref{thm:ltb}, item 2, when $\Theta$ is bounded] 


We consider a subset of times where $\rho$ fulfils a suitable lower bound, namely:
\begin{align*}
      \mathcal{T}&:=\acc{t\in[1,+\infty)\
   \tq \sup_{x\in \Er^N}\rho(t,x)\goq\frac{\lambda[r]}2}.
  \end{align*}
  Of course, it may happen that $\mathcal{T}=\emptyset$.
  We first derive a lower bound for $u$ on the set of times 
  $$\mathcal{T}+\{1\}:=\{t+1 \tq t\in\mathcal{T}\},$$
  using the parabolic Harnack inequality; next 
  we prove a lower bound for $u$ for times outside 
  $$\mathcal{T}+[0,1]:=\{t+s\tq t\in\mathcal{T},\ s\in[0,1]\},$$
  exploiting the principal eigenfunction; we finally 
  derive a lower bound for $\rho$ on $\mathcal{T}+[0,1]$.

  \paragraph{Step 1. Lower bound for $u$ for on $\mathcal{T}+\{1\}$.}

  By the definition of  $\mathcal{T}$ and the positivity of $\lambda[r]$,
  for any $t \in \mathcal{T}$
  there exists $x_t\in\Er^N$ such that $\rho(t,x_t) \geq \lambda[r]/4$.
  Letting $|\Theta|$ denote the measure of $\Theta$, one then deduces
  \[\forall t \in \mathcal{T},\quad
  \sup_{\theta\in\Theta}\,u(t,x_t,\theta)\goq
  \frac{\rho(t,x_t)}{|\Theta|}
  \goq\frac{\lambda[r]}{4|\Theta|}.\]
  Therefore, by the Harnack inequality \eqref{Harnack-bdd} derived before in the
    case $\Theta$ bounded, one infers that, for any $R>0$, there exists $C_R>0$ such that
  \begin{equation}\label{infu>}
  \forall t \in \mathcal{T},\quad
  \inf_{B_{\Er^N}(x_t,R)\times\Theta}u(t+1,\cdot,\cdot)\geq
  \frac1{C_R}\,\sup_{\Theta}\,u(t,x_t,\cdot)
  \goq\frac{\lambda[r]}{4C_R|\Theta|}.
  \end{equation}

  \paragraph{Step 2. Lower bound for $u$ on
  $[1,+\infty)\setminus(\mathcal{T}+[0,1])$.}
  First of all, we derive an estimate similar to \eqref{infu>} at time $1$.
  Namely, it follows from the parabolic strong maximum principle that 
  $u(t,x,\theta)>0$ for $t>0$, $x\in\Er^N$, $\theta\in\Theta$.
  Moreover, owing to the Hopf lemma, the Neumann boundary condition prevents 
  $u$ from vanishing on $\partial\Theta$. One deduces  in particular that
  $u(1,x,\theta)>0$ for $x\in\Er^N$, $\theta\in\overline\Theta$.
  This, together with \eqref{infu>}, implies that, for any $R>0$, there
  exists $C_R'>0$ such that
  \begin{equation}\label{eq:T+1}
  \forall t \in \mathcal{T}\cup\{0\},\quad
  \inf_{B_{\Er^N}(x_t,R)\times\Theta}u(t+1,\cdot,\cdot)\geq
  C_R',
  \end{equation}
  where we can choose, for instance, $x_0=0$.
  Next, by Proposition \ref{ppn:br_bounded} (and with the same notations), there exists $\bar R>0$,
  such that $\lambda_{\bar R}^m>\lambda[r]/2$.
  It follows that the associated $\varphi^m_{\bar R}$,
  extended to $0$ for $x\notin B_{\Er^N}(0,\bar R)$,
  is a generalised subsolution of the equation
  \begin{equation}\label{eq:persistence_linear}
    \dr_tv=\Delta_xv+\Delta_\theta v+v\pth{r(x,\theta)-\lambda[r]/2},
  \end{equation}
  and fulfils the Neumann boundary condition for $\theta\in\partial\Theta$.
  By periodicity, the same is true for any translation $\varphi^m_{\bar R}(\cdot-h,\cdot)$
  with $h\in\Zed^N$. 
  We normalise $\varphi^m_{\bar R}$ by $\|\varphi^m_{\bar R}\|_{\infty}=1$.
  
  Consider now a time $t\in[1,+\infty)\setminus(\mathcal{T}+[0,1])$, that is, 
  $t\geq1$ and $(t-1,t)\cap\mathcal{T}=\emptyset$.
  Call
  $$t_0:=\max\{s\leq t-1 \tq s\in\mathcal{T}\cup\{0\}\}.$$
  On the one hand, we know from \eqref{eq:T+1} that 
  \begin{equation*}
  \inf_{B_{\Er^N}(x_{t_0},\bar R+\sqrt N)\times\Theta}u(t_0+1,\cdot,\cdot)\goq C_{\bar R+\sqrt N}'.
  \end{equation*}
  On the other hand, taking $h\in\Zed^N$ (depending on $t_0$) such that 
  $x_{t_0}-h\in\mcc=[0,1]^N$, we have that the function 
  $\varphi^m_{\bar R}(\cdot-h,\cdot)$ vanishes outside 
  $B_{\Er^N}(h,\bar R)\subset
  B_{\Er^N}(x_{t_0},\bar R+\sqrt N)$ and therefore
  $C_{\bar R+\sqrt N}'\,\varphi^m_{\bar R}(\cdot-h,\cdot)$ lies below $u$ at time $t_0+1$.
  Finally, since $(t_0+1,t)\cap \mathcal{T}=\emptyset$, it follows from
  the definition of $\mathcal{T}$ that $u$ is a subsolution to~\eqref{eq:persistence_linear}
  on $(t_0+1,t)\times\Er^N\times\Theta$.
  We can then apply the parabolic comparison principle and infer that 
  $C_{\bar R+\sqrt N}'\,\varphi^m_{\bar R}(\cdot-h,\cdot)$ lies below $u$ for all times
  in $[t_0+1,t]$, hence in particular
  \[
  \sup_{x\in\Er^N}\Big(\min_{\theta\in\overline{\Theta}}u(t,x,\theta)\Big)
  \goq C_{\bar R+\sqrt N}'\,
  \min_{\theta\in\overline{\Theta}}\varphi^m_{\bar R}(0,\theta)>0.\] 
   This immediately gives a uniform 
    lower bound for $\sup_{x\in\Er^N}\rho(\cdot,x)$ on
    $[1,+\infty)\setminus(\mathcal{T}+[0,1])$.
    Thus, in order to conclude the proof, it remains to
    get an estimate on $\mathcal{T}+[0,1]$.

\paragraph{Step 3. Lower bound for $\rho$ on $\mathcal{T}+[0,1]$.}
Integrating with respect to $\theta$ the weak formulation~\eqref{eq:weak_formulation},
one ends up with the equation \eqref{eq:rho_pde} for $\rho$. 
We rewrite it~as
 \[
  \dr_t\rho=\Delta_x\rho+\tilde r( t,x)\rho-\rho^2,
    \qquad t>0,\ x\in\Er^N,
   \]
  where \[\tilde r(t,x):=\frac{1}{\rho(t,x)}\int_{\Theta}r(x,\theta)u(t,x,\theta)\de\theta\] satisfies $\inf r\leq \tilde r\leq \sup r$ (recall that 
  $r$ is bounded in the case $\Theta$ bounded).
%
 Standard parabolic estimates  and Morrey's inequality then yield
 $\rho\in{C}^{0,\alpha}([1,+\infty)\times\Er^N)$ 
   (for any $\alpha\in(0,1)$). Hence $\rho$ is uniformly continuous on $[1,+\infty)\times\Er^N$.
     
Consider $t\in \mathcal{T}$ and, as in the Step 1, let $x_t$ be such that
$\rho(t,x_t)>\lambda[r]/4$.
  By the uniform continuity of $\rho$, there exists $\delta>0$ independent of $t$
such that, for all $t\in\mct$, we have $\rho(t,\cdot)>\lambda[r]/8$
on $B_{\Er^N}(x_t,\delta)$. Therefore, calling $w$ the solution to
  $$\dr_t w=\Delta_x w+(\inf r)w-w^2,
    \qquad t>0,\ x\in\Er^N,$$
with initial datum equal to $\lambda[r]/8$ on $B_{\Er^N}(0,\delta)$
and $0$ outside,
we have that 
$$\min_{s\in[t,t+1]}\rho(s,x_t)\geq \min_{s\in[0,1]}w(s,0),$$
which is positive due to the strong maximum principle,
and independent of $t$.
This is the uniform lower bound on $\mathcal{T}+[0,1]$ which concludes the proof.
\end{proof}

\begin{proof}[Proof of Theorem \ref{thm:ltb}, item 2, when $\Theta=\Er^P$]
   Owing to Lemma \ref{lem:rho<},
   we are able to prove the result following essentially the same strategy 
   as in the case $\Theta$ bounded.
   However, several modifications are needed. 

 First of all, by the positivity of $\lambda[r]$ and 
hypothesis~\eqref{eq:r<0}, there exists $M>0$
such that 
\begin{equation}\label{r<outside}
r(x,\theta)<\frac{\lambda[r]^2}{8\|\rho\|_\infty} \qquad\text{for all }\;
x\in\Er^N,\ \|\theta\|\goq M
\end{equation}
(recall that $\rho$ is bounded, as well as $u$).
We call as usual $\bar r:=\sup r$. We then define 
  \begin{align*}
      \mathcal{T}&:=\acc{t\in[1,+\infty)\
   \tq \sup_{x\in \Er^N}
   \int_{\|\sigma\|\loq M}u(t,x,\sigma)\de \sigma\goq
   \frac{\lambda[r]^2}{ 8\bar r}}.
  \end{align*}
  As before, but in a slightly more general setting,
  we will derive a uniform lower bound for $u$ on the set of times 
  $\mathcal{T}+\{T\}$, depending on $T>0$. 
  In contrast with the
  case $\Theta$ bounded, using a principal eigenfunction of the linearised operator,
  we will be able to obtain a lower bound for $u$ outside
  $\mathcal{T}+[0,T]$ only for
  $T$ sufficiently large.
  We will then conclude by an estimate for~$u$ on $\mathcal{T}+[0,T]$.

  \paragraph{Step 1. Lower bound for $u$ on $\mathcal{T}+\{T\}$,
  for $T>0$.}
  For $t \in \mathcal{T}$, we take $x_t\in\Er^N$ such that 
  $$\int_{\|\sigma\|\loq M}u(t,x_t,\sigma)\de \sigma\goq
   \frac{\lambda[r]^2}{ 16\bar r}.$$
  We use the standard interior Harnack inequality which gives, for any
  $T,R>0$, a constant $C_{T,R}>0$ such that
  \[
  \forall t \in \mathcal{T},\quad
  \sup_{ B_{\Er^P}(0,R)}\,u(t,x_t,\cdot)\loq
  C_{T,R}\,\inf_{B_{\Er^{N+P}}((x_t,0),R)}u(t+T,\cdot,\cdot).\]
  One then finds, for $R\geq M$, 
  \begin{align*}
  \forall t \in \mathcal{T},\quad
  \inf_{B_{\Er^{N+P}}((x_t,0),R)}u(t+T,\cdot) &\goq
  \frac1{C_{T,R}}\;\sup_{ B_{\Er^P}(0,R)}\,u(t,x_t,\cdot)\\
  &\goq
  \frac1{C_{T,R}}\;\sup_{ B_{\Er^P}(0,M)}\,u(t,x_t,\cdot)\\
  &\goq\frac1{C_{T,R}|B_{\Er^P}(0,M)|}\int_{\|\sigma\|\loq M}u(t,x_t,\sigma)\de \sigma\\
  &\goq\frac{\lambda[r]^2}{ 16\bar r C_{T,R}|B_{\Er^P}(0,M)|}\,.
  \end{align*}

\paragraph{Step 2. Lower bound for $u$ on
  $[T,+\infty)\setminus(\mathcal{T}+[0,T])$, for $T$
  sufficiently large.} 
  Consider an arbitrary quantity $T>0$, that we will fix later. 
  We start by deriving a lower bound for~$u$
  at time $T$.
  This is simply given by the parabolic strong maximum principle,
  which ensures that 
  $u>0$ on $(0,+\infty)\times\Er^N\times\Er^P$.
  Thus, owing to Step 1, we deduce that for any $T>0$ and $R\geq M$, there
  exists $C_{T,R}'>0$ such~that
  \begin{equation}\label{eq:T+1whole}
  \forall t \in \mathcal{T}\cup\{0\},\quad
  \inf_{B_{\Er^{N+P}}((x_t,0),R)}u(t+T,\cdot,\cdot)\geq
  C_{T,R}',
  \end{equation}
  where we have set, say, $x_0:=0$. 
  
  Next, by Proposition \ref{ppn:cv_approximation_lambda} there exists $\bar R\geq M$
  such that $\lambda_{\bar R}^d>3\lambda[r]/4$.
  It follows that the associated $\varphi^d_{\bar R}$,
  extended to $0$ outside $B_{\Er^{N+P}}(0,\bar R)$,
  is a generalised subsolution of the equation 
  \begin{equation}\label{eq:persistence_linear_unbdd}
    \dr_tv=\Delta_xv+\Delta_\theta v+v\pth{r(x,\theta)-3\lambda[r]/4}.
  \end{equation}
  By periodicity, the same is true for each of its translations $\varphi^d_{\bar R}(\cdot-h,\cdot)$ by $h\in\Zed^N$. 
  We normalise~$\varphi^d_{\bar R}$ by $\|\varphi^d_{\bar R}\|_{\infty}=1$.
  
 Consider now a time $t\in[T,+\infty)\setminus(\mathcal{T}+[0,T])$, that is, 
  $t\geq T$ and $(t-T,t]\cap\mathcal{T}=\emptyset$.
  Call
  $$t_0:=\max\{s\leq t-T \tq s\in\mathcal{T}\cup\{0\}\}.$$
  On the one hand, we know from \eqref{eq:T+1whole} that
  \begin{equation*}
  \inf_{B_{\Er^{N+P}}((x_{t_0},0),\bar R+\sqrt N)}u(t_0+T,\cdot,\cdot)\goq C_{T,\bar R+\sqrt N}'.
  \end{equation*}
  On the other hand, taking $h\in\Zed^N$ (depending on $t_0$) such that 
  $x_{t_0}-h\in\mcc=[0,1]^N$, we have that the function 
  $\varphi^d_{\bar R}(\cdot-h,\cdot)$ vanishes outside 
  $B_{\Er^{N+P}}((h,0),\bar R)\subset
  B_{\Er^{N+P}}((x_{t_0},0),\bar R+\sqrt N)$. It follows that
  $C_{T,\bar R+\sqrt N}'\,\varphi^d_{\bar R}(\cdot-h,\cdot)$ lies below $u$ at time $t_0+T$.
  
  Finally, in order to compare $C_{T,\bar R+\sqrt N}'\,\varphi^d_{\bar R}(\cdot-h,\cdot)$
  with $u$ at time $t$, we now show that if $T>0$ is chosen large enough, then $u$ is a supersolution to~\eqref{eq:persistence_linear_unbdd}
  on $[t_0+T,t]\times\Er^N\times\Er^P$.
  Take any $\tau\in[T,t-t_0]$. 
  We apply estimate~\eqref{rho<H}
  from Lemma~\ref{lem:rho<}. We recall that the expression at the 
  right-hand side there
  is increasing with respect to $H$.
In the present case, $H$ is
  bounded from above by
  $\lambda[r]^2/4$,
  thanks to \eqref{r<outside} and the definition of $\mathcal{T}$
  (we have indeed $(t_0,t_0+\tau)\subset (t_0,t)$, which does not intersect $\mathcal{T}$).
  This shows that
    \[
    \sup_{x\in\Er^N}\,\rho(t_0+\tau,x)\loq 
    \frac{\lambda[r]}2\,\coth\left(\frac{\lambda[r]}2\,\tau+
    \arcoth\pth{\frac{2\|\rho\|_\infty+\lambda[r]}{\lambda[r]}}\right).
    \]
   Choosing $T$ large enough so that $\coth(\lambda[r]T/2)\leq 3/2$,
   the above estimate yields
    \[\forall \tau\in[T,t-t_0],\qquad \rho(t_0+\tau,\cdot)\leq 3\lambda[r]/4,\]
    which implies in turn
    that $u$ is a supersolution to~\eqref{eq:persistence_linear_unbdd}
    on $[t_0+T,t]\times\Er^N\times\Er^P$.
    
  We can then apply the parabolic comparison principle and infer that 
  $C_{T,\bar R+\sqrt N}'\,\varphi^d_{\bar R}(\cdot-h,\cdot)$ lies below $u$ 
  on $[t_0+T,t]\times\Er^N\times\Er^P$, whence in particular
  \[
  \sup_{x\in\Er^N}\Big(\inf_{\theta\in B_{\Er^P}(0,R/2)}u(t,x,\theta)\Big)
  \goq C_{\bar R+\sqrt N}'\,
  \inf_{\theta \in B_{\Er^P}(0,R/2)}\varphi^d_{\bar R}(0,\theta)>0.\] 
   This is the uniform 
    lower bound for $\sup_{x\in\Er^N}\rho(\cdot,x)$ on
    $[T,+\infty)\setminus(\mathcal{T}+[0,T])$.

\paragraph{Step 3. Lower bound for $u$ on $\mathcal{T}+[0,T]$.}
Consider $t\in \mathcal{T}$. Let $x_t$ be as in the Step 1 and let
$\theta_t\in\Er^P$ be such that $\|\theta_t\|\leq M$ and
$$u(t,x_t,\theta_t)=
\max_{\|\theta\|\loq M}u(t,x_t,\theta).$$
Then one has
$$u(t,x_t,\theta_t)\geq
\frac1{|B_{\Er^P}(0,M)|}\int_{\|\sigma\|\loq M}u(t,x_t,\sigma)\de \sigma\goq
   \frac{\lambda[r]^2}{8\bar r|B_{\Er^P}(0,M)|}.$$ 
        In the equation in \eqref{eq:main},
          the term $\rho$ is globally bounded 
          and $r(x,\theta)$ is bounded uniformly in $x\in\Er^N$ and locally uniformly in
$\theta\in\Er^P$, so, due to interior parabolic estimates,
the function $u$ is uniformly continuous
on $[1,+\infty)\times\Er^N\times B_{\Er^P}(0,M)$.
Therefore,
there exists $\delta>0$ independent of $t$ 
such that 
\[\inf_{B_{\Er^{N+P}}((x_t,\theta_t),\delta)}u(t,\cdot,\cdot)\geq
\frac{\lambda[r]^2}{16\bar r|B_{\Er^P}(0,M)|}.\]
(We point out that the centre of the ball
$(x_t,\theta_t)$ does depend on $t$.)
In order to translate this lower bound into an estimate which holds on the time interval $[t,t+T]$ and which is independent of $t$, 
 we need to be careful due to the possible unboundedness
of $r$.
We circumvent this difficulty by considering the solution $w$
to the Dirichlet problem
  \[
  \begin{cases}
  \dr_t w=\Delta_x w+w\left(\displaystyle\inf_{\Er^N\times B_{\Er^{P}}(0,M+\delta)}r -\|\rho\|_\infty\right)
    & \ \text{ in } (0,+\infty)\times B_{\Er^{N+P}}(0,\delta),\\
    w=0 & \ \text{ over } (0,+\infty)\times \partial B_{\Er^{N+P}}(0,\delta),
    \end{cases}\]
with initial datum equal to $\frac{\lambda[r]^2}{16\bar r|B_{\Er^P}(0,M)|}$. 
Observe that for $(x,\theta)\in B_{\Er^{N+P}}(0,\delta)$ it holds that
$(x+x_t,\theta+\theta_t)\in \Er^N\times B_{\Er^{P}}(0,M+\delta)$, hence
the function $u(\cdot+t,\cdot+x_t,\cdot+\theta_t)$ is a supersolution of the above problem, and moreover it lies above $w$ at time $0$.
It follows that 
$$\min_{[t,t+T]\times B_{\Er^{N+P}}(0,\delta/2)}u(\cdot,\cdot+x_t,\cdot+\theta_t)\geq 
\min_{[0,T]\times B_{\Er^{N+P}}(0,\delta/2)}w,$$
for any $T>0$, and in particular for the value of $T$ provided by the Step 2.
By the strong maximum principle, the right-hand side is a positive number.
This is the estimate on the set of times
$\mathcal{T}+[0,T]$ concluding the proof.
\end{proof}

\section{Dependence of the principal eigenvalue on the spatial period and 
on the mobility}
We recall that $r^L$ stands for the $L$-periodic version of $r$, i.e.,
\[r^L(x,\theta):=r\pth{\frac{x}{L},\theta}.\]

\begin{proof}[Proof of Proposition \ref{ppn:deriv_lambda}]
  Let us start  with the statement $(i)$ of the proposition.
  We first consider the case $\Theta=\Er^P$.
  For $L>0$, we let $\mcc^L:=[0,L]^N$ denote the periodic cell of size~$L$.
  For $R>0$, we call for short $B_R:=B_{\Er^P}(0,R)$ the ball of $\Er^P$ with radius $R$ and centre $0$.
  We let $(\varphi^L_R,\lambda_R^L)$ 
  satisfy the periodic eigenproblem analogous to~\eqref{eq:periodic_pevr}, but with $L$-periodicity instead of $1$-periodicity, that is, 
  \begin{equation*}
  \begin{cases}
    \mcl[r^L]\varphi^L_R=\lambda_R^L\varphi^L_R & \text{in }\Er^N\times B_R,\\
    \varphi^L_R=0 & \text{over $\Er^N\times\partial B_R$},\\
    \varphi^L_R>0 & \text{in $\Er^N\times B_R$},\\
    \varphi^L_R&\text{is $L$-periodic in $x$},
  \end{cases}
  \end{equation*}
  with the normalisation $\dabs{\varphi^L_R}_{L^2(\mcc^L\times B_R)}=1$.
  Take $L,L'>0$. We make now use of the Rayleigh formula for~$\lambda^{L'}_R$
  in order to estimate it in terms of~$\lambda^{L}_R$. 
  The Rayleigh formula,
  which is classical in the present case because the eigenfunctions
  act on a compact set (due to the $x$-periodicity), implies that
  \[
    \lambda_R^{L'}\goq\cro{-\int_{\mcc^{L'}\times B_R}\abs{\nabla \phi}^2+\int_{\mcc^{L'}\times B_R}\phi^2r^{L'} }\pth{\int_{\mcc^{L'}\times B_R}\phi^2}^{-1},
  \]
  for any function $\phi\in 
  W_{2}^{1}(\mcc^{L'}\times B_R)\cap C^0(\Er^N\times \overline{B_R})$
  which is $L'$-periodic in $x$ and vanishes on $\Er^N\times\partial B_R$.
  We call $h(x,\theta):=(\frac L{L'}x,\theta)$ and we take 
  $\phi=\varphi^L_R\circ h$ (which is $L'$-periodic) in the formula.
  This gives:
  \begin{align*}
    \lambda_R^{L'}&\goq\cro{-\int_{\mcc^{L'}\times B_R}\pth{\pth{\frac{L}{L'}}^2\abs{\nabla_x\varphi^L_R}^2\circ h+\abs{\nabla_\theta\varphi^L_R}^2\circ h}+\int_{\mcc^{L'}\times B_R} (\varphi^L_R\circ h)^2r^{L'}}\\
    &\esp\esp\times\pth{\int_{\mcc^{L'}\times B_R}(\varphi^L_R\circ h)^2}^{-1}.
  \end{align*}
  By the change of variables $(x',\theta')=h\pth{x,\theta}$, we obtain, since $\dabs{\varphi^L_R}_{L^2(\mcc^L\times B_R)}=1$:
  \begin{align*}
    \lambda_R^{L'}&\goq-\int_{\mcc^L\times B_R}\pth{\pth{\frac{L}{L'}}^2\abs{\nabla_x\varphi^L_R}^2+\abs{\nabla_\theta\varphi^L_R}^2}+\int_{\mcc^L\times B_R}(\varphi^L_R)^2r^L\\
    &=\lambda_R^L+\pth{1-\pth{\frac{L}{L'}}^2}\int_{\mcc^L\times B_R}\abs{\nabla_x\varphi^L_R}^2.
  \end{align*}
  Hence, for $L'>L$, we have $\lambda_R^{L'}> \lambda^L_R$, that is,
  $L\mapsto\lambda_R^{L}$ is increasing for all $R>0$.
  Since by Proposition~\ref{ppn:cv_approximation_lambda} $(i)$ (which holds
  true with $1$-periodicity replaced by $L$-periodicity) 
  we have that $\lambda_R^{L}\to\lambda[r^{L}]$
  as $R\to+\infty$,
  it follows that $L\mapsto \lambda[r^L]$ is nondecreasing.

  Let us now prove that $L\mapsto \lambda[r^L]$ is continuous. 
   We have shown before that 
\[
    \lambda_R^L-\lambda_R^{L'}\loq \pth{\pth{\frac{L}{L'}}^2-1}\int_{\mcc^L\times B_R}\abs{\nabla_x\varphi^L_R}^2,\]
therefore it will suffice to derive a bound for 
$\|\nabla_x\varphi^L_R\|_{L^2(\mcc^L\times B_R)}$
which is locally uniform with respect to $L>0$ and uniform
with respect to $R\geq 1$.
This follows by integrating on $\mcc^{L}\times B_R$ 
the equation satisfied by $\varphi^L_R$ multiplied
by $\varphi^L_R$ itself. Namely, recalling that 
$\dabs{\varphi^L_R}_{L^2(\mcc^L\times B_R)}=1$, one gets
  \begin{align*}
    \int_{\mcc^{L}\times B_R}\abs{\nabla_x\varphi^{L}_R}^2&=\int_{\mcc^{L}\times B_R}(r^{L}-\lambda_R^{L})(\varphi^{L}_R)^2
    \loq \big(\sup r-\lambda_R^{L}\big).
  \end{align*}
  Then, since $\lambda_R^{L}$ is increasing with respect to $R$ by 
  Proposition \ref{ppn:cv_approximation_lambda} $(i)$, and
  it is nondecreasing with respect to $L$, as we showed before, we find,
  for any given $L_0>0$,
  \[\forall R\geq1,\ \forall L, L'\geq L_0, \qquad
    \lambda_R^L-\lambda_R^{L'}\loq \Big|\pth{\frac{L}{L'}}^2-1\Big|
    \big(\sup r-\lambda_1^{L_0}\big).
    \]
    Passing to the limit as $R\to+\infty$ we derive the same inequality
    for the difference $\lambda[r^L]-\lambda[r^{L'}]$,
    which yields the continuity of $L\mapsto \lambda[r^L]$.

   Let us consider now the case where $\Theta$ is a smooth bounded subset of $\Er^P$.
  The generalised principal eigenvalue $\lambda[r^L]$ coincides now with
  the principal eigenvalue of
  the Neumann problem on $\Er^N\times\Theta$, 
  under $L$-periodicity condition in~$x$, which is then 
  a standard compact setting.
  We let $\varphi^L$ denote the (unique up to a scalar multiple)
  corresponding principal eigenfunction, provided by 
  Proposition \ref{ppn:br_bounded} $(ii)$.  
    Again, the Rayleigh formula
  is classical in such a framework. It implies 
  \begin{align*}
    \lambda[r^{L'}]&\geq\cro{-\int_{\mcc^{L'}\times \Theta}\abs{\nabla \phi}^2+\int_{\mcc^{L'}\times \Theta}\phi^2r^{L'} }\pth{\int_{\mcc^{L'}\times \Theta}\phi^2}^{-1},
  \end{align*} 
  for any function $\phi\in 
  W_{2}^{1}(\mcc^{L'}\times \Theta)$
  which is $L'$-periodic in $x$.
  Taking $\phi=\varphi^L\circ h$, with the same function $h$ as above, 
  one infers in the same way
  that $L\mapsto\lambda[r^L]$ is continuous and nonincreasing.

  Therefore, item $(i)$ is proved in both cases $\Theta$ bounded or $\Theta=\Er^P$.

  \medskip
    Let us now prove item $(ii)$. For a $x$-periodic function $\phi(x,\theta)$, we let $\phi^L(x,\theta):=\phi(x/L,\theta)$ be its $L$-periodic version. 
    We have
    \[\pth{ \mcl[r^L]\phi^L }(Lx,\theta)=\pth{ \mcl_{L^{-2}}[r] }\phi(x,\theta).\]
    It then follows from Definition~\ref{dfi:eigenvalue} that
    the eigenvalue $\lambda_d[r]$ associated with the operator $\mcl_{d}[r]$
    satisfies
    \[\lambda_d[r]=\lambda[r^{d^{-1/2}}].\]
    Thus, item $(ii)$ is a consequence of item $(i)$.

\end{proof}

\section*{Acknowledgements} 
This work has received support from the ANR project ReaCh, {ANR-23-CE40-0023-01}.
N.B. was supported by the Chaire Modélisation Mathématique et Biodiversité (École Polytechnique, Muséum national d’Histoire naturelle, Fondation de l’École Polytechnique, VEOLIA Environnement).
L.R. was supported by the 
European Union -- Next Generation EU, on the PRIN project
2022W58BJ5
 ``PDEs and optimal control methods in mean field games, population dynamics and multi-agent models''
 and by INdAM--GNAMPA.

\printbibliography

\end{document}